\documentclass{article}
\usepackage[utf8]{inputenc}
\usepackage{natbib}
\usepackage{graphicx}
\usepackage{amsmath}
\usepackage{amsfonts}
\usepackage{algorithm}
\usepackage{color}
\usepackage{algpseudocode}
\usepackage{comment}
\usepackage{verbatim}
\usepackage[most]{tcolorbox}
\title{Principled Graph Management}
\author{Julian Yarkony\\ Laminar Optimization Research Group, La Jolla, CA\\
Department of Computer Science, University of California, Irvine
\\ \\ Amelia Regan \\
Department of Computer Science, University of California, Irvine
}
\date{February 2022}

\usepackage{natbib}
\usepackage{graphicx}

\begin{document}

\maketitle

\begin{abstract}
   
   Graph Generation is a recently introduced enhanced Column Generation algorithm for solving expanded Linear Programming relaxations of mixed integer linear programs without weakening the expanded relaxations which characterize these methods. To apply Graph Generation we must be able to map any given column generated during pricing to a small directed acyclic graph for which any path from source to sink describes a feasible column. This structure is easily satisfied for vehicle routing, crew scheduling and various logistics problems where pricing is a constrained shortest path problem. Such graphs are then added to the restricted master problem when the corresponding column is generated during pricing. The construction of such graphs trades off the size/diversity of a subset of columns modeled by the graphs versus the additional computational time required to solve the restricted master problem induced by larger graphs.

    
    Graph Generation (GG) has two computational bottlenecks. The first is pricing. Pricing in GG and Column Generation (CG) is identical because of the structure of the problems solved. The second bottleneck is the restricted master problem (RMP), which is more computationally intensive in GG than in CG given the same number of columns generated. By design GG converges in fewer iterations than CG, and hence requires fewer calls to pricing. Therefore when the computation time of GG is dominated by pricing, as opposed to solving the RMP, GG converges much faster than CG in terms of time. However GG need not converge faster than CG when the GG RMP, rather than pricing, dominates computation.
    
    In this paper we introduce Principled Graph Management (PGM), which is an algorithm to solve the GG RMP rapidly by exploiting its special structure. We demonstrate the effectiveness of PGM inside a GG solution to the classical Capacitated Vehicle Routing Problem. Specifically we consider a regime where the number of  customers is large and heuristic pricing is used so that solving the GG RMP requires far more time than pricing. We demonstrate that PGM solves the GG RMP \textbf{hundreds} of times faster than the baseline solver and that the improvement in speed increases with problem size. 

\end{abstract}
\section{Introduction}
Expanded linear programming (LP) relaxations provide much tighter relaxations for broad classes of mixed integer linear  programs (MILP) found in logistics \citep{Desrochers1992,costa2019}, and more recently computer vision/machine learning \citep{yarkony2020data,branchimportant,HPlanarCC,FlexDOIArticle} than corresponding compact LP relaxations. Thus the solution of expanded LP relaxations provide the need for fewer branching operations in a branch-bound tree relative to compact LP relaxations.  Furthermore the fractional solution to an expanded LP relaxation can often be projected to a high quality feasible integer solution \citep{FlexDOIArticle}; which is not easily achieved using compact LP relaxations \citep{PlanarCC}. However solving expanded LP relaxations is non-trivial as the number of variables (referred to as columns) grows exponentially in the size of the original compact LP relaxation. In fact the set of columns in extended formulations is often not easily enumerated, much less able to be considered in optimization. Solving expanded LP relaxations is attacked using Column Generation (CG) \citep{barnprice, Desaulniers2005, gilmore1965multistage, cuttingstock}, which imitates the revised simplex approach \citep{evans1973revised}. CG solves expanded LP relaxations by generating a small but sufficient subset of the columns such that this subset is guaranteed to provide an optimal solution to the expanded LP relaxation. CG operates by alternating between solving the expanded LP relaxation over a limited subset of the columns (the restricted master problem or RMP) followed by adding in one or more negative reduced cost columns computed in a step called pricing. Pricing is typically a combinatorial optimization problem that is often a resource constrained shortest path problem (RCSPP), which can be attacked by special developed solvers \citep{Desaulniers2005}. CG terminates when no negative reduced cost columns exist, at which point the LP solution generated in the most recent RMP is optimal. Typically CG generates a very small portion of the columns that could be examined.  

CG is known to perform well in problems that fall under this regime:
\begin{itemize} 
\item We are given a set of agents that must cover a set of tasks. 
\item Each agent completes a subset of the tasks called an assignment. 
\item The problem is to provide each agent with an assignment so as to ensure that each task is covered at least once, and the total cost of the assignments is minimized. 

\end{itemize}
Note that the cost/structure of feasible assignments is not described here. CG can be used to solve these problems if we can compute the lowest reduced cost assignment (column) given rewards (negative costs) associated with completing each task during pricing.

The expanded LP relaxation is tighter than a compact relaxation when the pricing problem produces fractional solutions (where the fractional solutions have lower reduced cost than any integer solution) when solved trivially as a linear program. This property is studied under the name of the integrality property \citep{geoffrion1974lagrangean,vanderbeck2000dantzig, desrosiers2005primer}. 

CG demonstrates slow convergence in problems where the number of non-zero entries in the constraint matrix associated with columns generated during pricing exceeds a threshold number (such as 8-12 tasks in an assignment), as discussed in \citep{elhallaoui2005dynamic}). This difficulty can be circumvented with strong dual stabilization \citep{du1999stabilized,ben2006dual,elhallaoui2005dynamic,haghani2020smooth}. CG is accelerated in various ways, often by limiting the dual search space. One class of approaches concerns the use of application specific dual optimal inequalities (DOI), which bind the dual solution to lie in an easily described and dramatically smaller space than the original dual space; that space provably includes at least one optimal solution to the original expanded LP relaxation. In the primal problem DOI correspond to the addition of slack variables that are provably inactive in an optimal solution to the expanded linear programming relaxation. The use of DOI does not alter the structure of the pricing problem \citep{ben2006dual,haghani2020smooth}. Thus algorithms used to solve the pricing problem when no DOI are used can also be employed with the DOI.   

This paper builds off of Graph Generation \citep{yarkony2021graph} (GG), which can be understood as a dual stabilization scheme for CG for problems where pricing is a RCSPP. GG differs from CG only in the mechanism used to the produce a primal/dual solution at each iteration of GG.  
In essence GG produces a primal/dual solution at each iteration of GG that considers a large numbers of additional columns that are easily encoded in the RMP without exploding the computational difficulty of the  RMP. GG achieves this by adding graphs to the RMP for which every path from source to sink corresponds a column. Each column generated during pricing is mapped to a graph by GG that is designed to include a diverse set of columns related in some manner to the one generated during pricing. This graph is then added to the GG RMP. Thus each column generated during pricing permits the GG RMP to describe a combinatorial number of columns efficiently. GG does not alter the structure of the pricing problem from that of CG. GG requires only a problem domain specific means to map a column produced during pricing to such a graph, which is easy to do in broad classes of optimization problems. The use of larger graphs permits the description of larger more diverse sets of columns in the GG RMP at the cost of greater RMP computation time.  

By design GG decreases the total number of iterations of pricing required to solve the expanded LP relaxation (relative to CG). However this improvement comes at the expense of each solution to the RMP being more computationally intense relative to CG given the same number of columns generated. In cases where pricing rather than the RMP is the computational bottleneck for GG then it converges much faster than CG in terms of time. This is demonstrated in \citep{yarkony2021graph} for the case of the Capacitated Vehicle Routing Problem (CVRP). However the GG need not perform well when the RMP rather than pricing is the core computational bottleneck.  

 GG can be understood as a generalization of the work in \citep{de2002lp}, which introduces a compact LP relaxation called an arc based formulation for bin-packing/cutting stock problems that is exactly as tight as the well studied expanded LP relaxation; and is much tighter than the standard compact LP relaxation for bin-packing/cutting stock problems. GG extends the ability of such arc based formulations to cover more general classes of problems by introducing graphs associated with subsets of columns.

In this paper we show that the GG RMP is of special structure permitting much faster solution via a specialized approach that we call Principled Graph Management (PGM). PGM solves the GG RMP to circumvent the need to consider massive graphs in the RMP. PGM constructs a sufficient variables (called edges, referring to edges in the GG graphs) in the RMP by iterating between: \textbf{(1)} solving the RMP over the graphs with only the current set of edges (which is called the PGM-RMP) \textbf{(2)} adding edges in graphs associated with columns of negative reduced cost to the current set of edges under consideration. The determination of the edges to add is made using a fast dynamic program that computes the shortest path from source to sink including any given edge jointly. Only edges on the lowest reduced cost column on a given graph are added to the PGM-RMP. Note that this operation is a shortest path problem on a directed acyclic graph but not a RCSPP. To hot-start PGM at each iteration of GG we initialize the set of edges in the PGM-RMP to those with strictly positive values in the optimal solution to the RMP solution in the previous iteration. This ensures that we start PGM from the previous solution objective without having a large set of edges under consideration. PGM is named after Column Management, which removes/adds columns from the RMP in CG so as to prevent the CG RMP from becoming too computationally expensive \citep{desrosiers2005primer}.  \\

    \textbf{In summary, PGM exploits the following properties of the GG RMP:} 
    \begin{enumerate}
        \item the  addition of a small number of edges to a GG graph has the effect of adding a large number of columns of the expanded LP relaxation to the GG RMP without significantly increasing GG RMP solution time; 
        \item most edges are not active in the final solution to the GG RMP at any iteration of GG; 
        \item computing the lowest reduced cost column associated with any graph is a simple shortest path problem and not a resource constrained shortest path problem as in pricing. PGM solves the GG RMP in a manner akin to CG by alternating between the following two steps: 
        \begin{enumerate} 
        \item solving the RMP over a subset of the edges in the graphs;  \item adding edges to the subset under consideration associated with the lowest reduced cost column in each graph, which can be computed as a simple dynamic program. PGM terminates when the optimal solution to the GG RMP is solved, which is achieved when no negative reduced cost columns exist in the graphs. 
        \end{enumerate}
        \end{enumerate}

This paper is organized in the following way. In Section \ref{sec_Lit_Rev} we review the primary literature related to accelerating column generation (CG). In Section \ref{sec_review_CG} we review expanded LP relaxations and their corresponding solutions via CG. In Section \ref{sec_review_GG} we review the new technique which we call Graph Generation (GG). In Section \ref{Sec_capVVGG} we review the Capacitated Vehicle Routing Problem (CVRP) and its solution via GG. In Section \ref{Sec_Opt} we present our Principled Graph Management (PGM) approach for solving the GG restricted master problem (RMP). In Section \ref{Sec_exper} we provide experiments demonstrating the applicability of  PGM to a GG solution for CVRP, and by extension to many large scale combinatorial optimization problems. In Section \ref{Sec_conc} we conclude and discuss both experimental and computational extensions which will permit solving industrial scale problems.   
\section{Literature Review}
\label{sec_Lit_Rev}
In this section we review related work on accelerating Column Generation (CG).  
\subsection{Trust Region Based Stabilization}
Trust region based methods exploit the understanding of CG operating as a search algorithm over the dual space \citep{marsten1975boxstep} seeking to maximize the Lagrangian bound of the master problem (MP) (where MP is an alternative name for expanded LP relaxation). The Lagrangian bound at a given point in dual space is the lower envelope of a set of affine linear functions evaluated at that point in dual space. Here each affine function corresponds to the set of columns generated during pricing at that point in dual space. Distant points in dual space are often associated with very different binding affine functions so the restricted master problem (RMP) is often inclined to travel to points in dual space where its current set of columns in the RMP do not well describe the Lagrangian bound at the new point. This difficulty is circumvented by enforcing or encouraging the dual solution generated to remain near the dual point with greatest Lagrangian relaxation thus far identified.  Trust region based methods can be understood a compromise between standard CG and sub-gradient based methods \citep{barahona1998plant}, which take small steps in dual space so as to attempt avoid travel decreasing the Lagrangian bound significantly.
%
\subsection{Dual Optimal Inequalities}
Dual optimal inequalities (DOI) \citep{ben2006dual} provide provable bounds on the space where the optimal dual solution lies, thus decreasing the space over which CG searches. These bounds are easily computed based on problem class specific structure, and problem instance specific information. For example Smooth-DOI exploit the observation that for problems embedded on a metric space that the dual values must change smoothly over that space  \citep{haghani2020smooth}. By enforcing smoothness of the dual variables at each iteration of CG, convergence to the optimal solution is obtained in many fewer iterations. The use of DOI does not alter the structure of the CG pricing problem.  

\subsection{Advanced Pricing}
Often pricing can generate multiple negative reduced cost columns that cover distinct subsets of the dual variables. The addition of some or all of these columns to the RMP at each iteration of CG may make CG converge in far fewer iterations. Sometimes this can be done by solving parallel pricing problems, which must produce distinct solutions as in biological cell segmentation \citep{zhang2017efficient}. In cases where dynamic programming is used to solve pricing, many distinct columns can be generated by looking at the dynamic programming tables and returning multiple columns covering distinct (but not necessarily non-overlapping) subsets of the dual variables. Such methods have to trade off the desire to add more columns to the RMP, thus decreasing the number of iterations, with the increase in time for solving the RMP with additional columns. Column management \citep{desrosiers2005primer} helps to keep the RMP computationally efficient by repeatedly removing columns that are of little value, such as those that have not been active in the RMP solution for many iterations.  

Heuristic pricing can be used to rapidly generate low reduced cost columns instead of relying on an expensive (slow) exact solver \citep{Desaulniers2005} (chapter 4). Heuristic pricing can also be replaced with exact pricing when heuristic pricing fails to generate a negative reduced cost column. Thus the use of heuristic pricing need not preclude exact solution of the master problem.

\section{Review of Column Generation}
\label{sec_review_CG}

In this section we review expanded linear programming (LP) relaxations, and the standard Column Generation (CG) solution \citep{barnprice,lubbecke2005selected} for such problems. We use $\Omega$ to denote the set of primal variables (in the expanded linear program), which we index by $l$. The primal variables in the expanded LP relaxation are referred to as columns as is common in the CG literature. We should note that $\Omega$ is commonly too large to enumerate much less consider in optimization. We use $A_{:l}$ to denote the vector in the constraint matrix $A$ corresponding to column $l$. We use $c_l \in \mathbb{R}$ to denote the cost associated with column $l$. Typically $c_l$ is an easily evaluated function of $A_{:l}$. We use $\vec{b}$ to denote the vector of right hand side constants on the constraints, with associated dual variables denoted $\pi$. We use $\theta_l$ to denote the non-negative decision variable associated with column $l$.  
We write the standard expanded LP relaxation below, which we refer to as the master problem (MP).  
\begin{subequations}
\label{orig_LP}
\begin{align}
    \min_{\theta \geq 0} \sum_{l \in \Omega }c_l\theta_l\\
    \sum_{l \in \Omega}A_{:l}\theta_l \geq \vec{b} \quad [\pi]
\end{align}
\end{subequations}
Solving \eqref{orig_LP} is done using CG when $\Omega$ can not be explicitly enumerated (as is the case in many applications in operations research \citep{barnprice,lubbecke2005selected}); and the lowest reduced cost column in $\Omega$ can be computed given any non-negative dual solution $\pi$. CG can be used when any negative reduced cost column can be computed, if a negative cost column exists. CG can be used to approximately solve \eqref{orig_LP} when a heuristic method for identifying a negative reduced cost column is used; where the heuristic method may not return a negative reduced cost column if such a column exists.   
The computation of the lowest reduced column in $\Omega$ is written as follows using $\bar{c}_l$ to denote the reduced cost of column $l$.
\begin{subequations}
\label{pricing}
\begin{align}
    \min_{l \in \Omega}\bar{c_l}\\
    \bar{c}_l=c_l-A^{\top}_{:l}\pi \quad \forall l \in \Omega
\end{align}
\end{subequations}
%
CG is an exact solver for \eqref{orig_LP} that does not explicitly consider all of $\Omega$; since considering all of $\Omega$ is intractable. CG solves \eqref{orig_LP} in an iterative manner that alternates between the following two steps: \textbf{(a)} solving \eqref{orig_LP} over a limited subset of $\Omega$ denoted $\Omega_R$; and \textbf{(b)} adding negative reduced cost columns to $\Omega_R$ via solving \eqref{pricing}, which is typically a combinatorial optimization problem such as a resource constrained shortest path problem \citep{baldacci2011new,costa2019}. These two steps are referred to as solving the restricted master problem (RMP), and solving pricing respectively. We refer to the RMP as $\Psi(\Omega_R)$, which we write in primal/dual form below.   

\begin{align}
\label{orig_RMP}
    \Psi(\Omega_R)=\min_{\substack{\theta \geq 0\\ \sum_{l \in \Omega_R}A_{:l}\theta_l \geq \vec{b}}}\quad \sum_{l \in \Omega_R }c_l\theta_l
    =\max_{\substack{\pi \geq 0\\ c_l-A_{:l}^{\top}\pi \geq 0 \quad \forall l \in \Omega_R}}b^{\top}\pi
\end{align}
We initialize $\Omega_R$ with a set of columns describing a feasible solution to the RMP.  Commonly $\Omega_R$ is initialized using artificial variables with prohibitively high cost. When no negative reduced cost columns exist, given the dual optimal RMP solution then we have provably provided the optimal solution to the Master Problem. 
%
%
In Alg \ref{basicCG} we provide the CG solution to \eqref{orig_LP} in pseudo-code form with annotation below.    
\begin{algorithm}[!b]
 \caption{Basic Column Generation}
\begin{algorithmic}[1] 
\State $\Omega_R\leftarrow $ from user \label{alg_1_get_omR}
\Repeat  \label{alg_1_start_loop}
\State  $\theta,\pi\leftarrow $Solve \eqref{orig_RMP} over $\Omega_R$ \label{alg_1_RMP}
\State $l_* \leftarrow \mbox{arg}\min_{l \in \Omega}\bar{c}_l$ \label{alg_1_gen_col}
\State $\Omega_R \leftarrow \Omega_R \cup l_*$ \label{alg_1_add_col}
 \Until{$\bar{c}_{l_*} \geq 0$} \label{alg_1_end_loop}
 \State Return last $\theta$  generated \label{alg_1_ret}
\end{algorithmic}
\label{basicCG}
\end{algorithm} 

\begin{itemize}
    \item Line \ref{alg_1_get_omR}:  We initialize CG with columns from the user, which may consist of artificial variables with prohibitively high cost that ensure a feasible solution exists to the RMP.  
    \item Lines \ref{alg_1_start_loop}-Line \ref{alg_1_end_loop}:  We generate a sufficient set of columns to solve the MP exactly. 
    \begin{itemize}
        \item Line \ref{alg_1_RMP}:  Produce the solution to the RMP in \eqref{orig_RMP}.
        \item Line \ref{alg_1_gen_col}:  Produce the lowest reduced cost column.  CG does not require that \eqref{pricing} is solved exactly.  As long as a negative reduced cost column is added to $\Omega_R$ during any round of pricing (if such a column exists) then CG is guaranteed to solve \eqref{orig_LP} exactly.  Often more than one column is added to $\Omega_R$ during pricing.  This is facilitated in cases where pricing is solved using a dynamic program as dynamic programs generate many solutions over the course of optimization \citep{costa2019,wang2017tracking}.
        \item Line \ref{alg_1_add_col}:  We add the generated column to $\Omega_R$.  If we generated more than one column we would add some or all of these columns.  
        \item Line \ref{alg_1_end_loop}:  If no negative reduced cost column exists we terminate optimization.
    \end{itemize}
    \item Line \ref{alg_1_ret}:  We return the optimal solution.  This can be provided as input to a branch-price solver which calls Alg \ref{basicCG} in the inner loop.  
\end{itemize}
\section{Review of Graph Generation}
\label{sec_review_GG}
In this section we review the Graph Generation (GG) algorithm \citep{yarkony2021graph}. GG is an enhanced Column Generation (CG) algorithm that differs from CG by solving a more computationally intensive restricted master problem (RMP) at each iteration.  
GG is designed to require far fewer iterations that CG to converge to the optimal solution. Thus GG is faster than CG when solving pricing not solving the RMP dominates the computation time for GG. GG does not alter the structure of the pricing problem; nor does it loosen the master problem (MP).

GG requires the notion of a family of columns, which we now describe. Let $F$ be a set of subsets of $\Omega$, each member of which is called a family of columns.  We index $F$ using $f$ where $\Omega_f \subseteq \Omega$ is the set of columns in family $f$. There is a surjection from $\Omega$ to $F$ where for any $l \in \Omega$ the term $f_l$, which lies in $F$, is the index of the family corresponding to $l$.  For any $l \in \Omega$ the associated family $\Omega_{f_l}$ includes $l$ (meaning $l \in \Omega_{f_l} \quad \forall l \in \Omega$).  We should note that for any $l \in \Omega$ that $l$ may lie in many families not just $f_l$.  Each family is associated with the following terms.  
\begin{itemize}
    \item Let directed acyclic graph $G^f$ be associated with edge set $E^f$ and vertex set $V^f$. We index $E^f$ with $i,j$ corresponding to vertices in $V^f$.  There are special vertices  $v^+$ and $v^-$ in $V^f$ called the source and sink respectively. 
    \item Let $P^f$ be the set of paths in $G^f$ starting at $v^+$ and ending at $v^-$, which we index by $p$.  There is a surjection from $P^f$ to $\Omega_f$, where $p$ maps to $l_p$.  We describe paths using $a_{ijp} \in \{0,1\}$ where  $a_{ijp}=1$ IFF path $p$ includes edge $ij \in E^f$ and otherwise set $a_{ijp}=0$.  
    \item We now express the convex cone of paths, which we refer to as path cone using non-negative values $\{ \psi^f_{ij} \quad \forall ij \in E^f\}$ where $\psi_{ij}^f$ is the number of times edge $ij$ is covered in this combination.  The path cone is the set of possible settings of vector $\psi^f$ (over all $ij \in E^f$) corresponding to a non-negative combination of paths.  We describe path cone mathematically below. 
    \begin{align}
    \sum_{ij\in E^f}\psi^f_{ij}-\sum_{ji \in E^f}\psi^f_{ij}=0 \quad \forall i \in V^f-(v^+ \cup v^-)
\end{align}
    \item Each edge in $ij \in E^f$; is equipped with a cost $c^f_{ij}$ s.t. for any $p \in P^f$ the total cost of edges on the path $p$ is identical to that of $c_{l_p}$.  We write this formally below.
    \begin{align}
        c_{l_p}=\sum_{ij \in E^f}c^f_{ij}a_{ijp} \quad \forall p \in P^f
    \end{align}
    Thus we can map a non-negative vector $\psi^f$ that lies in the path cone to the corresponding contribution to the objective of the RMP.  
    \item Each edge $ij \in E^f$ is equipped with a vector $\vec{h}^f_{ij}$ with number of rows equal to the number of rows of $A$ s.t. the following holds.
    \begin{align}
        A_{:l_p}=\sum_{ij \in E^f}\vec{h}^f_{ij}a_{ijp} \quad \forall p \in P^f
    \end{align}
    Thus we can map a non-negative vector $\psi^f$ that lies in the path cone to the corresponding contribution to the constraints of the RMP.  
    \end{itemize}
Using the notation of graphs associated with families of columns we now consider the RMP that is solved at each iteration of GG.  Given any $\Omega_R \subseteq \Omega$, let $\Omega_{R2} $ be the union of columns in the families of columns in $\Omega_R$; meaning $\Omega_{R2} = \cup_{l \in \Omega_R}\Omega_{f_l}$. With the aim of accelerating the convergence of optimization GG solves $\Psi(\Omega_{R2})$ instead of over $\Psi(\Omega_{R})$  as in CG, at each iteration of GG.  This helps solve the MP faster (in terms of time) when solving this $\Psi(\Omega_{R2})$ leads to a smaller number of iterations of pricing being required to solve the MP, and $\Psi(\Omega_{R2})$ can be solved efficiently. Below we write the efficient solution to $\Psi(\Omega_{R2})$ using $F_R=\cup_{l \in \Omega_R}f_{l}$ where $F_R$ is the set of families associated with the columns generated thus far.
\begin{subequations}
\label{aug_RMP3}
\begin{align}
\Psi(\Omega_{R2})=    \min_{\substack{\psi \geq 0 }}\sum_{\substack{f \in F_R\\ij \in E^f}}c^f_{ij} \psi^f_{ij}\\
    \sum_{\substack{f \in F_R\\ij \in E^f}}\vec{h}^f_{ij}\psi^f_{ij} \geq \vec{b} \quad [\pi]\\
    \sum_{\substack{j\in V^f\\  ij \in E^f}} \psi^f_{ij}=\sum_{\substack{j\in V^f\\  ji \in E^f}} \psi^f_{ji} \quad \forall i \in V^f-(v^+,v^-), f \in F_R
\end{align}
\end{subequations}
The primal optimization formulation in \citep{yarkony2021graph} which proposed GG uses the $\theta$ terms from \eqref{orig_RMP} in addition to the $\psi$ terms. However the use of $\theta$ terms is redundant since $\Omega_{f_l}$ contains $l$ by construction.  We should note that \eqref{aug_RMP3} is a primal block angular linear program  where there is a bijection from $F_R$ to a block of variables.  Hence \eqref{aug_RMP3} can be attacked using by principled methods exploiting this structure \citep{castro2007interior}.  In Alg \ref{advCG} we describe GG in pseudo-code form, and provide annotation below.
\begin{algorithm}[!b]
 \caption{Graph Generation Algorithm (GG)}
\begin{algorithmic}[1] 
\State $F_R \leftarrow $ from user \label{alg_2_init_1} \label{alg_2_init_2}
\Repeat \label{alg_3_loop_start}
\State   $\psi,\pi \leftarrow$ Solve $ \Psi(\Omega_{R2})$ via \eqref{aug_RMP3}  \label{alg_2_solve_RMP}
\State $l_* \leftarrow \mbox{arg}\min_{l \in \Omega}\bar{c}_l$ \label{alg_2_pricing}
\State $F_R \leftarrow F_R \cup f_{l_*}$ \label{alg_2_add_F}
 \Until{$\bar{c}_{l_*} \geq 0$} \label{alg_3_loop_end}
 \State Return last $\psi$  generated.  \label{alg_3_return}
\end{algorithmic}
\label{advCG}
\end{algorithm} 

\begin{itemize}
    \item Line \ref{alg_2_init_1}: We initialize GG with one more more families from the user that together can describe a feasible solution RMP.  This may consist of families containing artificial columns with prohibitively high cost.  
    \item Line \ref{alg_3_loop_start}-\ref{alg_3_loop_end}:  We generate a sufficient set of graphs to solve the MP exactly.  
    \begin{itemize}
        \item Line \ref{alg_2_solve_RMP}:  Solve the RMP providing a primal/dual solution.  This RMP is called the GG RMP.  
        \item Line \ref{alg_2_pricing}:  Call pricing to generate the lowest reduced cost column $l_*$.  As in basic CG we may choose to generate more than one column.
    \item Line \ref{alg_2_add_F}: Add the new family $f_{l_*}$ associated with the column generated during pricing to the RMP.  If multiple columns are generated during pricing then a family is added to $F_R$ for each such column.
    \item Line \ref{alg_3_loop_end}:  If no column has negative reduced cost then we terminate optimization.      \end{itemize}

    \item Line \ref{alg_3_return}: Return the last primal solution generated, which is optimal.  As in standard CG this can be returned to a branch-price algorithm, that uses Alg \ref{advCG} as an inner loop operation.   
\end{itemize}
\section{Review of Graph Generation for Capacitated Vehicle Routing}
\label{Sec_capVVGG}

In this section we describe the application Graph Generation (GG) as it is used for the Capacitated Vehicle Routing Problem (CVRP) in \citep{yarkony2021graph}.  In CVRP we are given  set of customers that must be serviced each with a location and demand; the location of a depot where vehicles start and end; and a set of homogeneous vehicles with fixed capacity. We seek to assign vehicles to routes so as to minimize the total distance traveled; such that all customers are serviced, and the capacity of the vehicles is respected. CVRP is often attacked using Column Generation (CG) \citep{costa2019,Desrochers1992}. Here the pricing problem is a resource constrained shortest path problem (RCSPP), which is NP-hard \citep{costa2019}, with exact computation time growing exponentially in the number of resources.  \citep{Desrochers1992}.  This RCSPP problem has one resource for capacity remaining and one for each customer.  

We now provide a formal description of CVRP and the standard solution via Column Generation (CG). We define the set of customers as $N$, which we index by $u$.  We use $N^+$ to denote $N$ augmented with the depot. We have access to $K$ homogeneous vehicles, which start (and end) at the depot, each with capacity $d_0 \in \mathbb{Z}_+$.  The demand of customer $u$ is denoted $d_u \in \mathbb{Z}_{+}$. We use $c_{uv}$ to denote the distance between any pair $u,v$ each of which lie in $N^+$.  We use $\Omega$ to denote the set of feasible routes (columns in the master problem (MP)), which we index by $l$. A route is feasible if it starts/ends at the depot, services no more demand than $d_0$, and visits each customer one or zero times (but never a fractional number or a number greater that one).  We describe $l$ using the following notation.  
\begin{itemize}
    \item We set $a_{ul}=1$ if route $l$ services customer $u$, and otherwise set $a_{ul}=0$ for any $u \in N$.  
    \item We set $a_{uvl}=1$ if $v$ is immediately preceded by $u$ in the route $l$, and otherwise set $a_{uvl}=0$ for any $u\in N^+,v \in N^+$. 
    \item 
    Below we define $c_l$ to be the cost of route $l$ where $c_l$ is defined as the total travel distance on route $l$ using $c_{uv}$, which is the distance from $u$ to $v$.
    \begin{align}
        c_l=\sum_{\substack{u \in N^+\\ v \in N^+}}a_{uvl}c_{uv} \quad \forall l \in \Omega
    \end{align}
\end{itemize}
We use $\theta_l$ to denote the decision variable for the number of times route $l$ is selected in our solution.  
We write the MP for CVRP formally below (with exposition provided below the equations) and dual variables $\pi_u,\pi_0$ written in $[]$ next to their associated equations. 
\begin{subequations}
\label{RMP_CVRPTW}
\begin{align}
\min_{\theta \geq 0}\sum_{l \in \Omega}c_l\theta_l \label{eq_rmp_cvrp_obj}\\
\sum_{l \in \Omega}a_{ul}\theta_l \geq 1 \quad \forall u \in N  \quad [\pi_u] \label{eq_rmp_cvrp_cover}\\
\sum_{l \in \Omega}-\theta_l\geq -K \quad [\pi_0]  \label{eq_rmp_cvrp_pack}
\end{align}
\end{subequations}

In \eqref{eq_rmp_cvrp_obj} we seek to minimize the total distance traveled over all vehicles used.  In \eqref{eq_rmp_cvrp_cover} we enforce that each customer is serviced at least once.  We should observe that in no optimal solution to the MP is a customer covered more than once. However in optimal solutions to the RMP this may not hold, and customers may be covered more than once.  In the dual form this corresponds to enforcing that \eqref{eq_rmp_cvrp_cover} is associated with a non-negative dual variable, which decreases the dual search space for CG, thus accelerating the convergence of CG.  In \eqref{eq_rmp_cvrp_pack} we enforce that no more than $K$ vehicles are used.

We now consider the application of GG to CVRP.  We associate each family $f$ with  an ordered list containing all $u \in N$ using $\beta^f_u \in \mathbb{Z}_{+}$ where $\beta^f_u$ is the position in the ordered list that $u$ occupies. A column $l \in \Omega$ lies in $\Omega_f$ if for any $u\in N,v\in N$ s.t.  $a_{uvl}=1$ then $u$ must come before $v$ in the ordered list associated with family $f$. 
%
Formally we write this as follows. 
\begin{align}
    \quad (l \in \Omega_f)\leftrightarrow \{ (a_{uvl}=1) \rightarrow (\beta^f_u <\beta^f_v) \quad \forall f \in F, l \in \Omega, u \in N,v \in N \}
\end{align}
We now describe the GG graph $G^f=(V^f,E^f)$ as described in \citep{yarkony2021graph} for CVRP as follows.  For every $u \in N$, $d \in [0,1,2...,d_0-d_u]$ create one vertex in $V^f$ denoted $(u,d)$.  We construct $E^f$, and the corresponding cost terms $c_{ij}^f$ as follows.  Connect vertex $(u,d)$ to vertex $(v,d-d_v)$ with an edge of cost $c_{uv}$ IFF $\beta^l_u<\beta^l_v$ (and $d-d_v\geq 0$). Connect every vertex $(u,d)$ to the sink $v^-$ with an edge of cost $c_{u,-2}$; where $c_{u,-2}$ is the distance from $u$ to the depot.  For every $u \in N$ connect the source vertex $v^+$ to vertex $(u,d_0-d_u)$ with an edge of cost $c_{-1,u}$;  where $c_{-1,u}$ is the distance from the depot to $u$. 

 We now consider the construction of $\{\vec{h}^f_{ij} \quad \forall ij \in E^f\}$.  We index the constraints in the MP for CVRP in \eqref{RMP_CVRPTW} with $u\in N$ and $0$ for the the constraint enforcing no more than $K$ routes used.  We denote the corresponding entries in $\vec{h}^f_{ij}$ as $\vec{h}^f_{ij;u} $ and $\vec{h}^f_{ij;0}$ respectively.  For any $\vec{h}^f_{ij}$ where $j$ corresponds to $(u,d)$ for some $(u,d)$ we define $\vec{h}^f_{ij;u}=1$.  All other entries of $\vec{h}^f_{ij}$ are zero for that $ij \in E^f$.  For $j=v^-$ we define $\vec{h}^f_{ij;0}=-1$. All other entries of $\vec{h}^f_{ij}$ are zero for that $ij \in E^f$.


We now describe the generation of the orderings associated routes so that, given $l$ returned by pricing, that $\Omega_{f_l}$ contains useful routes to improve the objective of the RMP\citep{yarkony2021graph}.  Observe that customers that are in similar physical locations should be in similar positions on the ordered list so that a route in $\Omega_{f_l}$ can visit any subset of customers in a group of nearby customers without leaving the area, and then coming back, or forgoing the visiting of other distant customers. Of course it is impossible to preserve all such nearby spatial relationships if the world is not one dimensional (and in CVRP it is generally two dimensional).  However we seek to construct the ordering so as to encourage these relationships to be satisfied especially when one member of the relationship is a customer used in column $l$ (where the set of customers covered in route $l$ is denoted $N_l$).  We initialize the ordering with the customers in $N_l$ sorted in order from first visited to last visited on route $l$.  
 Now iterate over the remaining customers ($N-N_l$) in a random order.  For each customer $u$ insert it immediately behind the customer in $N_l$ that is closest to $u$.  For a customers closer to the depot than any customer in $N_l$ insert them in the beginning of the list.

We have already used quite a bit of notation in the preceding sections. To improve the readability of the most important section of our paper, we provide a table of some key notation here:

\begin{tcolorbox}

\textbf{NOTATION}\\

\textbf{From Section \ref{sec_review_GG}: Graph Generation:}\\


$\Omega$ is the set of feasible routes (columns) indexed by $l$ 

$F$ is a set of subsets of $\Omega$, each member of which is called a family, indexed by $f$ 

$\Omega_f \subseteq \Omega$ is the set of columns in family $f$ 

For any $l \in \Omega$ the associated family $\Omega_{f_l}$ includes $l$ 

$G^f$ is a directed acyclic graph with edge set $E^f$ and vertex set $V^f$ 

$E^f$ is indexed with $i,j$ corresponding to vertices in $V^f$ 

$v^+$ and $v^-$ in $V^f$ are the source and sink respectively   

$P^f$ is the set of paths in $G^f$ starting at $v^+$ and ending at $v^-$ indexed by $p$.  

$\{ \psi^f_{ij} \quad \forall ij \in E^f\}$ is a path cone, where $\psi_{ij}^f$ is the number of times edge $ij$ is covered in this combination.  

The path cone is the set of possible settings of vector $\psi^f$ (over all $ij \in E^f$) corresponding to a non-negative combination of paths.  

Given any $\Omega_R \subseteq \Omega$, $\Omega_{R2}$ is the union of columns in the families of columns in $\Omega_R$ 


\end{tcolorbox}

\section{Principled Graph Management}
\label{Sec_Opt}
Solving the Graph Generation (GG) restricted master problem (RMP) can become computationally  difficult as the number of families in the RMP grows over the course of GG iterations, and for problems where larger graphs are used. In \citep{yarkony2021graph} (which introduced GG) computation time for GG is so heavily dominated by pricing that GG RMP solution time of \eqref{aug_RMP3} is unimportant.  However in this document we consider problems where heuristic pricing is employed and the size of the problems is very large (and hence large graphs are desired) and thus computation time is dominated by solving the GG RMP (in \eqref{aug_RMP3}).  In this section we consider the fast solution to \eqref{aug_RMP3}, using an algorithm which we refer to as Principled Graph Management (PGM).  PGM constructs small subsets of the edges in the graphs corresponding to families in $F_R$ so that the RMP over those partial graphs provides the same solution as the complete RMP.  PGM achieves this in a manner akin to how CG solves the master problem (MP). 

We describe these partial graphs a follows. Let $\hat{G}^{f}=\hat{V}^f,\hat{E}^f,$ where $\hat{V}^f \subseteq V^f$ and $\hat{E}^f\subseteq E^f$ for each $f\in F_R$. Here $\hat{V}^f$ is defined to include each $j \in V^f$ s.t. there exists an edge including $j$ in $\hat{E}^f$.  The family of columns corresponding to paths in $\hat{G}^f$ is denoted $\hat{\Omega}_f$.  The union of the columns in such families is denoted $\hat{\Omega}_{R2}$ and defined formally as follows $\hat{\Omega}_{R2}=\cup_{f \in F_{R}} \hat{\Omega}_f$.  We write $\Psi(\hat{\Omega}_{R2})$ below as an LP.

\begin{subequations}
\label{aug_RMP_Reduced}
\begin{align}
\Psi(\hat{\Omega}_{R2})=    \min_{\substack{ \psi \geq 0 }}\sum_{\substack{f\in F_R\\ ij \in \hat{E}^f}}c^f_{ij} \psi^f_{ij}\\
    \sum_{\substack{f \in F_R\\ij \in \hat{E}^f}}\vec{h}^f_{ij}\psi^f_{ij} \geq b \quad [\pi] \\
    \sum_{\substack{j\in \hat{V}^f\\  ij \in 
    \hat{E}^f}} \psi^f_{ij}=\sum_{\substack{j\in \hat{V}^f\\  ji \in \hat{E}^f}} \psi^f_{ij} \quad \forall i \in \hat{V}^f-(v^+,v^-), f \in F_R
\end{align}
\end{subequations}
Given any dual solution to \eqref{aug_RMP_Reduced} let $\mu_{ij;f}$ denote the reduced cost of the lowest reduced cost column in $\Omega_f$ for which an associated path including edge $i,j$.  We define $\mu_{ij;f}$ explicitly below.
\begin{align}
    \mu_{ij;f}=\min_{\substack{p \in \mathcal{P}^f\\ a_{ijp}=1}}\bar{c}_{l_p} \quad  \forall ij \in E^f,f \in F_R
\end{align}
%

PGM generates sufficient sets $\{\hat{E}^f \quad \forall f \in F_R\}$ s.t. $\Psi(\hat{\Omega}_{R2})=\Psi(\Omega_{R2})$ by alternating between the following two steps in a manner akin to CG:  \textbf{(1)} Solving optimization in \eqref{aug_RMP_Reduced} using an off the shelf LP solver or by exploiting its primal block angular structure by using the solver of \citep{castro2007interior}; \textbf{(2)} Adding edges associated with negative reduced cost columns in $\Omega_{f} $ for each $f \in F_R$. Specifically we add to $\hat{E}^f$ any edges associated with the lowest reduced cost column (that has negative reduced cost) in family $f$.  Note that the addition of a small number of edges in $\hat{E}^f$ can vastly increase the number of columns in $\hat{\Omega}_f$ since every path in $E^f$ corresponds to a column in $\Omega_f$.  
So as to hot start PGM, we  initialize $\{\hat{E}^f \}\quad \forall f \in R $ to include all edges with positive value $(\psi^f_{ij}>0)$ in the most recent primal optimal solution to the GG RMP (meaning the result of Line \ref{alg_2_solve_RMP} of Alg \ref{advCG}).  We terminate PGM when no $l \in \cup_{f \in F_{R}}\Omega_{f}$ has negative reduced cost meaning $\mu_{ij;f}\geq 0$ for all $ij \in E^f,f \in F_R$ thus certifying that we have solved the RMP optimally.  In Alg \ref{PGM_basic_alg} we provide pseudo-code for PGM with exposition provided below.  
\begin{algorithm}[!b]
 \caption{Principled Graph Management(PGM)}
\begin{algorithmic}[1] 
\State $F_R \leftarrow $ from user \label{init_PGM1}
\State $\hat{E}^f\leftarrow $ from user \label{alg_get_E} for all $f \in F_R$ \label{init_PGM2}
\Repeat \label{alg_2_loop_start} \label{ittAlg}
\State   $\psi,\pi \leftarrow$ Solve $ \Psi(\hat{\Omega}_{R2})$ via \eqref{aug_RMP_Reduced} \label{alg_3_solve_RMP} 
\For{$f \in F_R$}\label{grabTermsStart}
\State Compute $\mu_{ij,f}$ via \eqref{my_eq_mu}  for all $ij \in E^f$\label{get_muEq}
\State $\hat{E}^f \leftarrow \hat{E}^f \cup \{ ij \in E^f; \mu_{ij,f}<0 $ and $\mu_{ij;f}=\min_{\hat{i}\hat{j}\in E^f}\mu_{\hat{i}\hat{j},f} \}$.   \label{add_edges}
\EndFor \label{grabTermsEnd}
 \Until{$\mu_{ij,f}\geq 0 \quad \forall ij \in E^f, f \in F_R$} \label{alg_2_loop_end}
 \State Return last $\psi,\pi$  generated.  \label{alg_2_return}
\end{algorithmic}
\label{PGM_basic_alg}
\end{algorithm} 

\begin{itemize}
    \item Line \ref{init_PGM1}-\ref{init_PGM2}:  Initialize $F_R$ and $\{\hat{E}^f, \forall f \in F_R\}$ to include a feasible solution to \eqref{aug_RMP3}.  The selected edges are those used in the optimal primal solution to the  previous solution to the RMP (if this is not the first iteration of GG).  If it is the first iteration of GG then any edges describing a feasible solution to the RMP are sufficient.  Such a choice can include artificial variables. 
    \item Line \ref{ittAlg}-\ref{alg_2_loop_end}:  Construct sufficient sets of edges $\hat{E}^f$ for all $f \in F_R$ s.t. $\Psi(\hat{\Omega}_{R2})=\Psi(\Omega_{R2})$
    \begin{itemize}
    \item Line \ref{alg_3_solve_RMP}: Solve $\Psi(\hat{\Omega}_{R2})$ providing a dual solution to the RMP over edge sets $\hat{E}^f$ for each $f \in F_R$
     \item Line \ref{grabTermsStart}-\ref{grabTermsEnd}: Add any edges  included in a column with lowest reduced cost (for that family); and that has negative reduced cost.
     \begin{itemize}
\item Line \ref{get_muEq}:  Compute $\mu_{ij;f}$ for each $i,j \in E^f$.  We show how to do this efficiently subsequently in this section.  
    \item Line \ref{add_edges}:  Add any edges to $\hat{E}^f$ for which $\mu_{ij;f}$ is minimized (and is negative valued).
     \end{itemize}
    \item Line \ref{alg_2_loop_end}:  Terminate the optimization when $\mu_{ij;f}\geq 0$ for all $ij \in E^f, f \in F_R$ meaning that we have produced an optimal solution to the RMP.  
    \end{itemize}
    \item Line \ref{alg_2_return}:  Return the solution of the RMP 

\end{itemize}
We now consider the fast computation of the $\mu_{ij;f}$ terms using the following helper terms.  We define  $\mu^+_{i,f}$  to be the cost of the lowest cost path starting at $v^+$ and ending at $i$ on directed acylcic graph $G^f$ with edge weights defined as follows  for each edge $ij \in E^f$:  $c^f_{ij}-(\vec{h}^f_{ij})^{\top}\pi$.  
Similarly we define $\mu^-_{i,f}$ to be the cost of the lowest cost path starting at $i$ and ending at $v^-$ on that graph. We now write $\mu_{ij,f}$ using recursive definitions of $\mu^+_{if},\mu^-_{jf}$ below.
\begin{subequations}
\label{my_eq_mu}
\begin{align}
        \mu_{ij,f}=\mu^+_{if}+\mu^-_{jf}+c^f_{ij}-(\vec{h}^f_{ij})^{\top}\pi \quad \forall ij \in E^f\label{mu_f_edge}\\
        \label{mu_f_+}
    \mu^+_{jf}=\min_{ij \in E^f}\mu^+_{if}+c^f_{ij}-(\vec{h}^f_{ij})^{\top}\pi \quad \forall j \in V^f-v^+\\
\label{mu_f_-}
    \mu^-_{if}=\min_{ij \in E^f}\mu^-_{jf}+c^f_{ij}-(\vec{h}^f_{ij})^{\top}\pi \quad \forall i \in V^f-v^-\\
    \mu^+_{v^+f}=\mu^-_{v^-f}=0
\end{align}
\end{subequations}
Since $E^f$ describes a directed-acyclic graph there exists an ordering (without ties) of the vertices that can be easily computed via expanding a vertex only after all of its predecessors are expanded.  
Thus we can compute \eqref{mu_f_+} in this order  in time in order $O(|E^f|)$. Similarly we can compute \eqref{mu_f_-} in the reverse of this order.  The values $\mu^+,\mu^-$ then can be used to produce the $\mu_{ij,f}$ terms.  Alternatively the solution to \eqref{mu_f_+} can be efficiently written as a shortest path calculation from $v^+$ to all vertices. Similarly the solution to \eqref{mu_f_-} can be efficiently written as a shortest path calculation from all vertices to $v^-$.

\section{Experiments}
\label{Sec_exper}
In this section we demonstrate that Principled Graph Management (PGM) accelerates the convergence of Graph Generation (GG).  To do this we consider a problem domain in which solving the restricted master problem (RMP) is the computational bottleneck for GG not pricing.  We then show that each iteration of RMP in GG is solved dramatically faster using PGM than solving the RMP with a baseline solver (denoted BL).  We refer to solving an expanded LP relaxation with GG using PGM/BL to solve the RMP as GG+PGM/GG+BL respectively.  We also show that GG+PGM is rapidly able to (approximately) solve the master problem (MP), while GG+BL is unable to do so efficiently. By approximately solve we mean that we have generated a primal-dual solution pair of identical objective for which the heuristic pricing solver is unable to generate a negative reduced cost column within a certain number of attempts.  We consider the problem domain of Capacitated Vehicle Routing Problem (CVRP) using a large number of customers.  In this domain the GG RMP is not efficiently solved by BL but is efficiently sovled by PGM.  We use the same mechanism to map columns to families as described in Section \ref{Sec_capVVGG}, and detailed in \citep{yarkony2021graph} 

We organize this section as follows. In Section \ref{sub_sec_data} we describe our problem instance data set.  In Section \ref{sub_sec_heur} we describe the heuristic pricing approach that we employ.  In Section \ref{sub_sec_init} we describe the initialization of the $F_R,\hat{E}^f$ for GG. In  \ref{sub_sec_lp} we describe our hardware/software used and implementation details.  
In Section \ref{sub_sec_res} we provide experimental results and analysis.

\subsection{Problem Instance Data Set}
\label{sub_sec_data}

We consider problem instances with $150$ customers placed randomly according uniform distribution on a $50$ by $50$ grid each of demand one. We have access to $40$ vehicles each of capacity $6$.  The depot is placed randomly as well according to a uniform distribution.  All distances between customers (and or the depot) use the $\ell_2$ distance which is rounded up to the nearest integer.  We generated ten such problem instances randomly.  

\subsection{Heuristic Pricing}
\label{sub_sec_heur}

To solve the CVRP pricing problem we use the heuristic described below based on the heuristic pricing scheme used in \citep{haghani2021multi}. This heuristic is parameterized by a random vector (describing a topological ordering of the customers) that restricts the set of columns (routes) that can be generated. Hence when our heuristic fails to generate a negative reduced cost column with a given random vector we try a different random vector.  We try up to $100$ random vectors before terminating pricing and declaring that GG has converged.  

The heuristic pricing approach is designed to avoid having to consider a resource constrained shortest path problem  (RCSPP) with many resources.  Exact pricing for CVRP is difficult because we must ensure that a route does not visit the same customer more than once. Our heuristic pricing scheme considers optimization over a subset of the routes such that this constraint need not be explicitly considered. Our heuristic is described as follows: First we generate a random topological ordering of the customers (that includes no ties).  We describe this ordering with $\kappa_{uv}=1$ if $u$ comes before $v$ (not necessarily immediately before) in this ordering and otherwise set $\kappa_{uv}=0$.  We generate the lowest reduced cost $l \in \Omega$ consistent with this ordering where consistency is defined as follows. For any pair of customers $u,v$ for which $a_{uvl}=1$ then it must be the case that $\kappa_{uv}=1$.  Solving for the lowest reduced cost $l\in \Omega$ consistent with this ordering is a resource constrained shortest path problem with one resource corresponding to the remaining capacity of the vehicle.  Pricing over this restricted set can be trivially written as a shortest path problem as described in \citep{haghani2021multi}.  We refer to this heuristic pricing approach as topological pricing since it is based on a random topological order of the customers.  

\subsection{Initialization of Graph Generation}
\label{sub_sec_init}

We now consider the initialization of $F_R$ and the associated set of edges.  We initialize $F_R$ with only one family of columns.  This family corresponds to the lowest reduced cost column generated if we weighted all customers equally with very high weight.  Specifically we take the solution to pricing as solved via heuristic pricing  where dual solution sets $\pi_u\leftarrow$ to a large positive number for all $u\in N$ and $\pi_0\leftarrow0$.  Pricing produces a low reduced cost column denoted $l$, which is associated with the family $f_l$.  Observe that $f_l$ can describe a feasible solution to the RMP.  This is because any subset of $d_0$ customers can be described in a column in $\Omega_{f_l}$ with some non-infinite cost corresponding those customers.  We then initialize $\hat{E}^{f_l}$ to be equal to $E^{f_l}$.  This is chosen because in the first iteration of GG the GG RMP is not computationally problematic in our experiments.  Other initialization choices could be made as we only need to ensure that a feasible solution to $\Psi(\hat{\Omega}_{R2})$ exists.   

\subsection{Hardware/Software Employed and Implementation Details}
\label{sub_sec_lp}

All experiments are conducted on a 2020 Dell computer with an Intel(R) Core(TM) i7-10750H CP \@2.60 giga hertz(GHZ) with 16 gigabytes of memory(GB).  All experiments are conducted using MATLAB version R2020 b (update 3).  All linear programs are solved with the MATLAB built in ``linprog" linear programming (LP) solver with default options, and no initial solution provided.  Computation of $\mu^+,\mu^-$ terms is done using the built in MATLAB ``distance" function implementing the Bellman-Ford algorithm \citep{ford2015flows}.  
We set the maximum computation time for each algorithm (GG+BL and GG+PGM) to solve the problem instance to be 3000 seconds.  Once this time limit has been reached we permit the given iteration of GG to be completed.  

To prevent numerical issues from occurring inside of edge addition for PGM we relax the edge addition criteria by adding a tiny offset to the edge selection rule for Line \ref{add_edges} of PGM (Alg \ref{PGM_basic_alg}). Thus we update  $\hat{E}^f$ as follows.   $\hat{E}^f \leftarrow \hat{E}^f \cup \{ ij \in E^f; \mu_{ij,f}<0 \quad \mbox{ and }\mu_{ij,f}< \min_{\hat{i}\hat{j} \in E^f}\mu_{\hat{i}\hat{j},f} +\epsilon\}$ for $\epsilon$ as a small positive number ($\epsilon=10^{-3}$).

\subsection{Experimental Results and Analysis}
\label{sub_sec_res}
We now compare the performance of GG+PGM vs GG+BL.  We compare the time required to solve each problem instance up till heuristic pricing fails to generate a negative reduced cost column (at which point the reduced cost of the column generated is regarded as zero and GG is terminated). 

\begin{figure}[!hbtp]
	\includegraphics[width=0.49\linewidth]{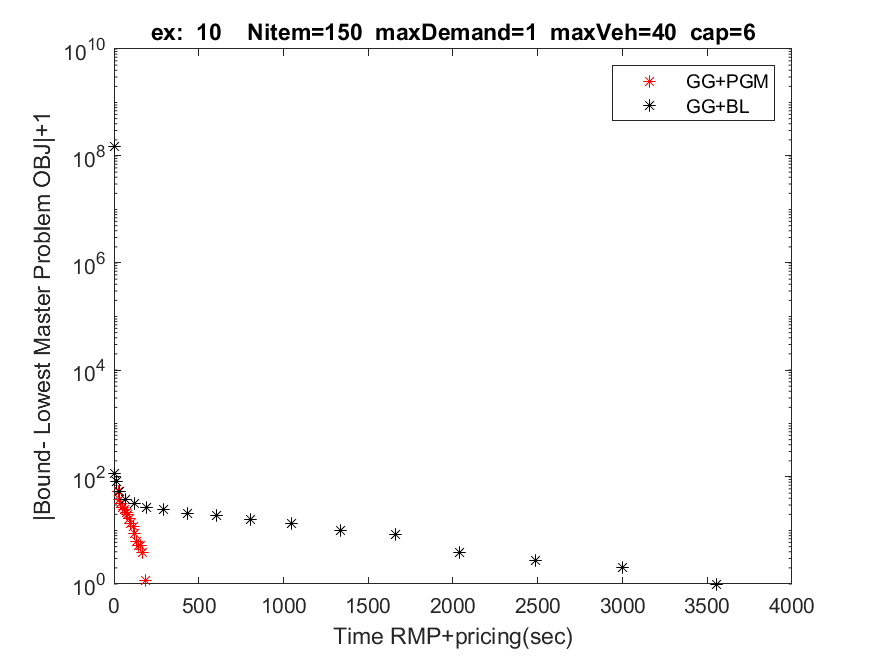}
	\includegraphics[width=0.49\linewidth]{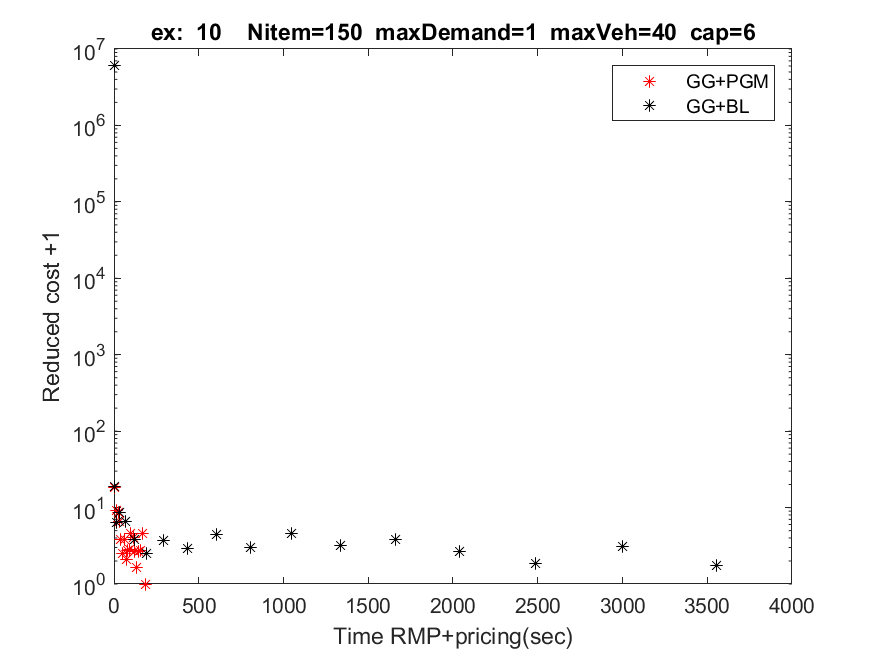}\\
	\caption{Results on an individual problem instance as a function of time (sec) in semi-log scale. Individual dots show the a value (LP RMP or -reduced cost) on $y$ axis for the given algorithm at the iteration for that dot.  The left side provides the results RMP objective and the right side minus 1 times the reduced cost of the column generated during pricing. We add one to the $y$ values of all terms which lets us use the semi-log scale.  
	}
	\label{fig:iter0}
\end{figure}

In Figure \ref{fig:iter0} we show the convergence on a specific example. In Figure \ref{fig:iter0} (left) we show the LP objective as a function of time for GG+PGM and GG+BL.  Each dot indicates the time that an iteration is completed, and the associated LP value. In Figure \ref{fig:iter0} (right)  we show the reduced cost of the column computed during pricing as a function of time for GG+PGM and GG+BL. Each dot indicates the time that an iteration is completed and the associated reduced cost. In Appendix \ref{additionalResults} we provide the performance for all remaining problem instances in our data set as they are provided in Fig \ref{fig:iter0}.

\begin{figure}[!hbtp]
	\includegraphics[width=0.49\linewidth]{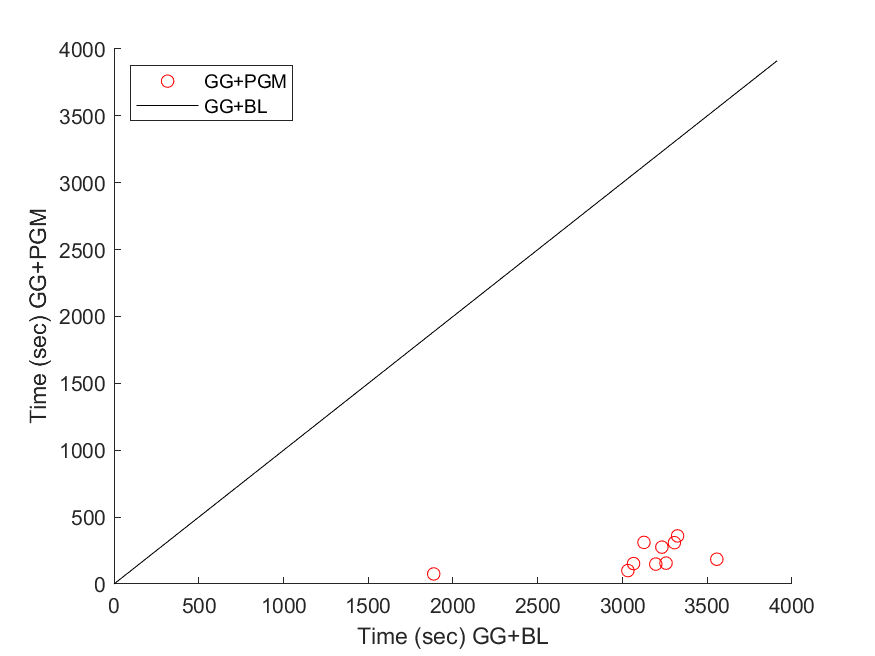}
		\includegraphics[width=0.49\linewidth]{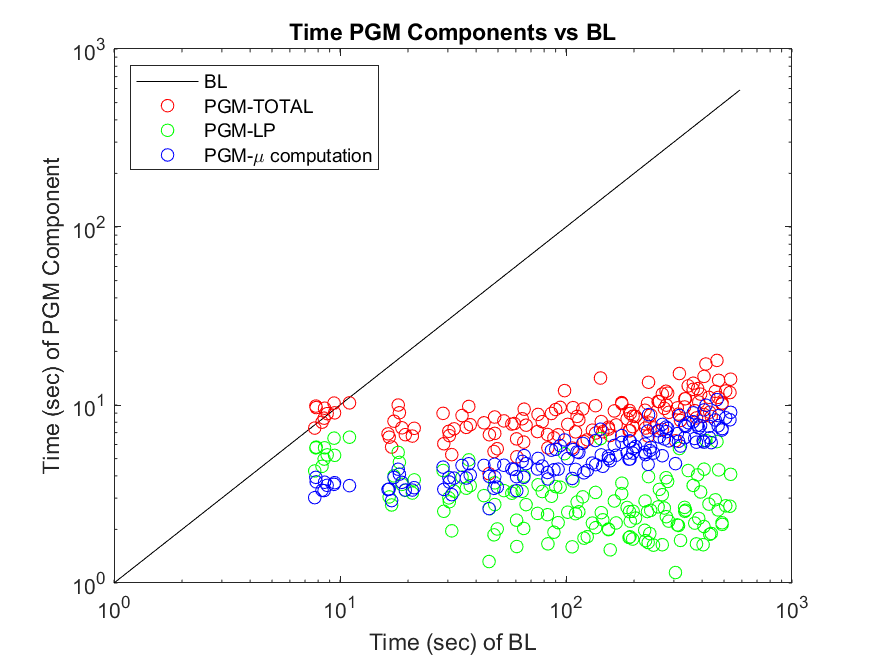}  
	\caption{Left:  Time till convergence of GG+BL vs GG+PGM over problem instances.  (Right): Time comparison for BL vs PGM.  Each dot describes the time to solve the RMP by BL vs PGM.  PGM time consumption is broken down into components.
	}
	\label{rmp_tiome}
\end{figure}

We plot the aggregated results of this experiment over all problem instances in our data set in Figure \ref{rmp_tiome} (left).  Each data point describes the time taken for that problem instance to be solved by GG+PGM vs GG+BL.  To provide a baseline for improvement, we plot a black line describing the linear function with slope =1 and $y$ intercept=0. The further below the line a dot is placed the greater corresponding performance improvement achieved by GG+PGM over GG+BL.  In Fig \ref{rmp_tiome} (left) observe that PGM converges faster than GG, and is able to (approximately) solve problem all problem instances within the time limit, while GG+BL is only able to approximately solve one problem instance.  

We also compared the time required to solve the RMP generated by GG+PGM at each iteration by PGM vs BL.  We provide these results in Figure \ref{rmp_tiome} (right) aggregated over problem instances and calls to solving the RMP. For each iteration of PGM/problem instance we provide three data points (in red,blue,green respectively); where the x coordinate corresponds to the time for the RMP to be solved using the BL.  The red dot's $y$ coordinate describes the total time for PGM to solve the RMP. The blue dot's $y$ coordinate describes the total time spent during PGM to solve the $\mu$ computation.  The green dot's $y$ coordinate describes the total time spent during PGM to solve the LPs. We observe that the time to solve the RMP differs widely for BL while being relatively constant and much smaller for PGM.  In nearly all cases PGM takes less time than BL; and in the cases that PGM take more time than BL are in the regime is the smallest time regime studied.  In the high computation time regime for GG we observe factor 50 speedups by using PGM over BL.  We also observe that time spent in PGM is dominated by the computation of the $\mu$ terms.  These can be efficiently computed exploiting parallel computation permitting us to only consider the green dots to measure performance.  In this case we observe speed factors up to $266$.  

\section{Conclusion and Future Research}
\label{Sec_conc}
In this paper we introduce a new mechanism to accelerate the solution of expanded linear programming (LP) relaxations via Graph Generation (GG), by dramatically speeding up the solution to the restricted master problem (RMP) of GG. We refer to our approach as Principled Graph Management (PGM). PGM avoids considering all variables in the GG RMP by constructing a sufficient set to exactly solve the GG RMP in a manner akin to traditional Column Generation (CG). PGM constructs a sufficient set of variables in the RMP by iterating between solving the RMP and by adding variables corresponding to edges associated with columns of lowest reduced cost on each graph.  The determination of the variables (edges) to add is a fast dynamic program, which computes the shortest path including any given edge jointly.

While we have incredibly promising results for both graph generation alone, and in combination with principled graph management, our next step is to test these methods on commonly used test problems and then on much larger industrial scale problems. While Capacitated Vehicle Routing Problems provide a valid testing ground for these methods, we intend to extend our exploration to other well known and well tested problems and problem instances.

Methodologically, in future work we seek to further accelerate optimization of the RMP by removing unnecessary edges from consideration during PGM. Specifically if we reach a point where the LP of PGM becomes too big to solve, we would remove edges associated with exclusively positive reduced cost columns.  
This should only be done after an iteration of GG+PGM improves the objective so as to ensure convergence.  
We also intend to investigate the use of PGM methods in the context of Detour-Dual Optimal Inequalities \citep{yarkony_Detour_DOI}. In this manner the set of possible detours considered in the RMP would be increased when such detours are associated with the lowest reduced cost column(s) with detours.  
We also intend to speed up the solution the RMP by exploiting the natural primal block angular structure of these problems using the approach of \citep{castro2007interior}. Our initial exploration and discussions with Professor Castro are very promising. 
For time window constrained problems we intend to implement the discretization techniques introduced in \citep{boland2017continuous} and used in \citep{haghani2020integer}.
\bibliographystyle{abbrvnat} 
\bibliography{col_gen_bib}

\begin{thebibliography}{33}
\providecommand{\natexlab}[1]{#1}
\providecommand{\url}[1]{\texttt{#1}}
\expandafter\ifx\csname urlstyle\endcsname\relax
  \providecommand{\doi}[1]{doi: #1}\else
  \providecommand{\doi}{doi: \begingroup \urlstyle{rm}\Url}\fi

\bibitem[Baldacci et~al.(2011)Baldacci, Mingozzi, and Roberti]{baldacci2011new}
R.~Baldacci, A.~Mingozzi, and R.~Roberti.
\newblock New route relaxation and pricing strategies for the vehicle routing
  problem.
\newblock \emph{Operations Research}, 59\penalty0 (5):\penalty0 1269--1283,
  2011.

\bibitem[Barahona and Jensen(1998)]{barahona1998plant}
F.~Barahona and D.~Jensen.
\newblock Plant location with minimum inventory.
\newblock \emph{Mathematical Programming}, 83\penalty0 (1):\penalty0 101--112,
  1998.

\bibitem[Barnhart et~al.(1996)Barnhart, Johnson, Nemhauser, Savelsbergh, and
  Vance]{barnprice}
C.~Barnhart, E.~L. Johnson, G.~L. Nemhauser, M.~W.~P. Savelsbergh, and P.~H.
  Vance.
\newblock Branch-and-price: Column generation for solving huge integer
  programs.
\newblock \emph{Operations Research}, 46:\penalty0 316--329, 1996.

\bibitem[Ben~Amor et~al.(2006)Ben~Amor, Desrosiers, and Val{\'e}rio~de
  Carvalho]{ben2006dual}
H.~Ben~Amor, J.~Desrosiers, and J.~M. Val{\'e}rio~de Carvalho.
\newblock Dual-optimal inequalities for stabilized column generation.
\newblock \emph{Operations Research}, 54\penalty0 (3):\penalty0 454--463, 2006.

\bibitem[Boland et~al.(2017)Boland, Hewitt, Marshall, and
  Savelsbergh]{boland2017continuous}
N.~Boland, M.~Hewitt, L.~Marshall, and M.~Savelsbergh.
\newblock The continuous-time service network design problem.
\newblock \emph{Operations Research}, 65\penalty0 (5):\penalty0 1303--1321,
  2017.

\bibitem[Castro(2007)]{castro2007interior}
J.~Castro.
\newblock An interior-point approach for primal block-angular problems.
\newblock \emph{Computational optimization and Applications}, 36\penalty0
  (2-3):\penalty0 195--219, 2007.

\bibitem[Costa et~al.(2019)Costa, Contardo, and Desaulniers]{costa2019}
L.~Costa, C.~Contardo, and G.~Desaulniers.
\newblock Exact branch-price-and-cut algorithms for vehicle routing.
\newblock \emph{Transportation Science}, 26(1), 2019.

\bibitem[De~Carvalho(2002)]{de2002lp}
J.~V. De~Carvalho.
\newblock L{P} models for bin packing and cutting stock problems.
\newblock \emph{European Journal of Operational Research}, 141\penalty0
  (2):\penalty0 253--273, 2002.

\bibitem[Desaulniers et~al.(2005)Desaulniers, Desrosiers, and
  Solomon]{Desaulniers2005}
G.~Desaulniers, J.~Desrosiers, and M.~M. Solomon, editors.
\newblock \emph{Column Generation}.
\newblock Springer, New York, 1st edition, 2005.

\bibitem[Desrochers et~al.(1992)Desrochers, Desrosiers, and
  Solomon]{Desrochers1992}
M.~Desrochers, J.~Desrosiers, and M.~Solomon.
\newblock A new optimization algorithm for the vehicle routing problem with
  time windows.
\newblock \emph{Operations Research}, 40\penalty0 (2):\penalty0 342--354, 1992.

\bibitem[Desrosiers and L{\"u}bbecke(2005)]{desrosiers2005primer}
J.~Desrosiers and M.~E. L{\"u}bbecke.
\newblock A primer in column generation.
\newblock In G.~Desaulniers, J.~Desrosiers, and M.~M. Solomon, editors,
  \emph{Column Generation}, pages 1--32. Springer, New York, NY, 2005.

\bibitem[Du~Merle et~al.(1999)Du~Merle, Villeneuve, Desrosiers, and
  Hansen]{du1999stabilized}
O.~Du~Merle, D.~Villeneuve, J.~Desrosiers, and P.~Hansen.
\newblock Stabilized column generation.
\newblock \emph{Discrete Mathematics}, 194\penalty0 (1-3):\penalty0 229--237,
  1999.

\bibitem[Elhallaoui et~al.(2005)Elhallaoui, Villeneuve, Soumis, and
  Desaulniers]{elhallaoui2005dynamic}
I.~Elhallaoui, D.~Villeneuve, F.~Soumis, and G.~Desaulniers.
\newblock Dynamic aggregation of set-partitioning constraints in column
  generation.
\newblock \emph{Operations Research}, 53\penalty0 (4):\penalty0 632--645, 2005.

\bibitem[Evans and Steuer(1973)]{evans1973revised}
J.~P. Evans and R.~E. Steuer.
\newblock A revised simplex method for linear multiple objective programs.
\newblock \emph{Mathematical Programming}, 5\penalty0 (1):\penalty0 54--72,
  1973.

\bibitem[Ford and Fulkerson(1962)]{ford2015flows}
L.~R. Ford and D.~R. Fulkerson.
\newblock \emph{Flows in Networks}.
\newblock Princeton university press, 1962.

\bibitem[Geoffrion(1974)]{geoffrion1974lagrangean}
A.~M. Geoffrion.
\newblock Lagrangean relaxation for integer programming.
\newblock In \emph{Approaches to integer programming}, pages 82--114. Springer,
  1974.

\bibitem[Gilmore and Gomory(1961)]{cuttingstock}
P.~Gilmore and R.~Gomory.
\newblock A linear programming approach to the cutting-stock problem.
\newblock \emph{Operations Research}, 9\penalty0 (6):\penalty0 849--859, 1961.

\bibitem[Gilmore and Gomory(1965)]{gilmore1965multistage}
P.~Gilmore and R.~E. Gomory.
\newblock Multistage cutting stock problems of two and more dimensions.
\newblock \emph{Operations Research}, 13\penalty0 (1):\penalty0 94--120, 1965.

\bibitem[Haghani et~al.(2020)Haghani, Li, Koenig, Kunapuli, Contardo, and
  Yarkony]{haghani2020integer}
N.~Haghani, J.~Li, S.~Koenig, G.~Kunapuli, C.~Contardo, and J.~Yarkony.
\newblock Integer programming for multi-robot planning: A column generation
  approach.
\newblock \emph{arXiv preprint arXiv:2006.04856}, 2020.

\bibitem[Haghani et~al.(2021{\natexlab{a}})Haghani, Contardo, and
  Yarkony]{haghani2020smooth}
N.~Haghani, C.~Contardo, and J.~Yarkony.
\newblock Smooth and flexible dual optimal inequalities.
\newblock \emph{Informs Journal on Optimization, in press, arXiv preprint
  arXiv:2001.02267}, 2021{\natexlab{a}}.

\bibitem[Haghani et~al.(2021{\natexlab{b}})Haghani, Li, Koenig, Kunapuli,
  Contardo, Regan, and Yarkony]{haghani2021multi}
N.~Haghani, J.~Li, S.~Koenig, G.~Kunapuli, C.~Contardo, A.~Regan, and
  J.~Yarkony.
\newblock Multi-robot routing with time windows: A column generation approach.
\newblock \emph{arXiv preprint arXiv:2103.08835}, 2021{\natexlab{b}}.

\bibitem[Leal-Taixe et~al.(2012)Leal-Taixe, Pons-Moll, and
  Rosenhahn]{branchimportant}
L.~Leal-Taixe, G.~Pons-Moll, and B.~Rosenhahn.
\newblock Branch-and-price global optimization for multi-view multi-target
  tracking.
\newblock In \emph{Proc. 25th Conference on Computer Vision and Pattern
  Recognition}, pages 1987--1994, Providence, Rhode Island, 2012.

\bibitem[Lokhande et~al.(2020)Lokhande, Wang, Singh, and
  Yarkony]{FlexDOIArticle}
V.~S. Lokhande, S.~Wang, M.~Singh, and J.~Yarkony.
\newblock Accelerating column generation via flexible dual optimal inequalities
  with application to entity resolution.
\newblock In \emph{Proceedings of the AAAI Conference on Artificial
  Intelligence}, volume~34, pages 1593--1602, 2020.

\bibitem[L{\"u}bbecke and Desrosiers(2005)]{lubbecke2005selected}
M.~E. L{\"u}bbecke and J.~Desrosiers.
\newblock Selected topics in column generation.
\newblock \emph{Operations Research}, 53\penalty0 (6):\penalty0 1007--1023,
  2005.

\bibitem[Marsten et~al.(1975)Marsten, Hogan, and
  Blankenship]{marsten1975boxstep}
R.~E. Marsten, W.~Hogan, and J.~W. Blankenship.
\newblock The boxstep method for large-scale optimization.
\newblock \emph{Operations Research}, 23\penalty0 (3):\penalty0 389--405, 1975.

\bibitem[Vanderbeck(2000)]{vanderbeck2000dantzig}
F.~Vanderbeck.
\newblock On dantzig-wolfe decomposition in integer programming and ways to
  perform branching in a branch-and-price algorithm.
\newblock \emph{Operations Research}, 48\penalty0 (1):\penalty0 111--128, 2000.

\bibitem[Wang et~al.(2017)Wang, Wolf, Fowlkes, and Yarkony]{wang2017tracking}
S.~Wang, S.~Wolf, C.~Fowlkes, and J.~Yarkony.
\newblock Tracking objects with higher order interactions via delayed column
  generation.
\newblock In \emph{Proc. 20th International Conference on Artificial
  Intelligence and Statistics}, pages 1132--1140, Fort Lauderdale, Florida,
  2017.

\bibitem[Yarkony and Fowlkes(2015)]{HPlanarCC}
J.~Yarkony and C.~Fowlkes.
\newblock Planar ultrametrics for image segmentation.
\newblock In \emph{Proc. 28th Advances in Neural Information Processing
  Systems}, pages 64--72, Montreal, Quebec, 2015.

\bibitem[Yarkony et~al.(2012)Yarkony, Ihler, and Fowlkes]{PlanarCC}
J.~Yarkony, A.~Ihler, and C.~Fowlkes.
\newblock Fast planar correlation clustering for image segmentation.
\newblock In \emph{Proc. 12th European Conference on Computer Vision}, pages
  1169--1176, Florence, Italy, 2012.

\bibitem[Yarkony et~al.(2020)Yarkony, Adulyasak, Singh, and
  Desaulniers]{yarkony2020data}
J.~Yarkony, Y.~Adulyasak, M.~Singh, and G.~Desaulniers.
\newblock Data association via set packing for computer vision applications.
\newblock \emph{Informs Journal on Optimization}, 2\penalty0 (3):\penalty0
  167--191, 2020.

\bibitem[Yarkony et~al.(2021{\natexlab{a}})Yarkony, Haghani, and
  Regan]{yarkony2021graph}
J.~Yarkony, N.~Haghani, and A.~Regan.
\newblock Graph generation: A new approach to solving expanded linear
  programming relaxations.
\newblock \emph{arXiv preprint arXiv:2110.01070}, 2021{\natexlab{a}}.

\bibitem[Yarkony et~al.(2021{\natexlab{b}})Yarkony, Haghani, and
  Regan]{yarkony_Detour_DOI}
J.~Yarkony, N.~Haghani, and A.~Regan.
\newblock Detour dual optimal inequalities for column generation with
  application to routing and location.
\newblock \emph{arXiv preprint}, 2021{\natexlab{b}}.

\bibitem[Zhang et~al.(2017)Zhang, Wang, Gonzalez-Ballester, and
  Yarkony]{zhang2017efficient}
C.~Zhang, S.~Wang, M.~A. Gonzalez-Ballester, and J.~Yarkony.
\newblock Efficient column generation for cell detection and segmentation.
\newblock \emph{arXiv preprint arXiv:1709.07337}, 2017.

\end{thebibliography}
\appendix
\section{Additional Results:  Description of the convergence of the Individual problem instances}
\label{additionalResults}
\begin{figure}[!hbtp]
	\includegraphics[width=0.49\linewidth]{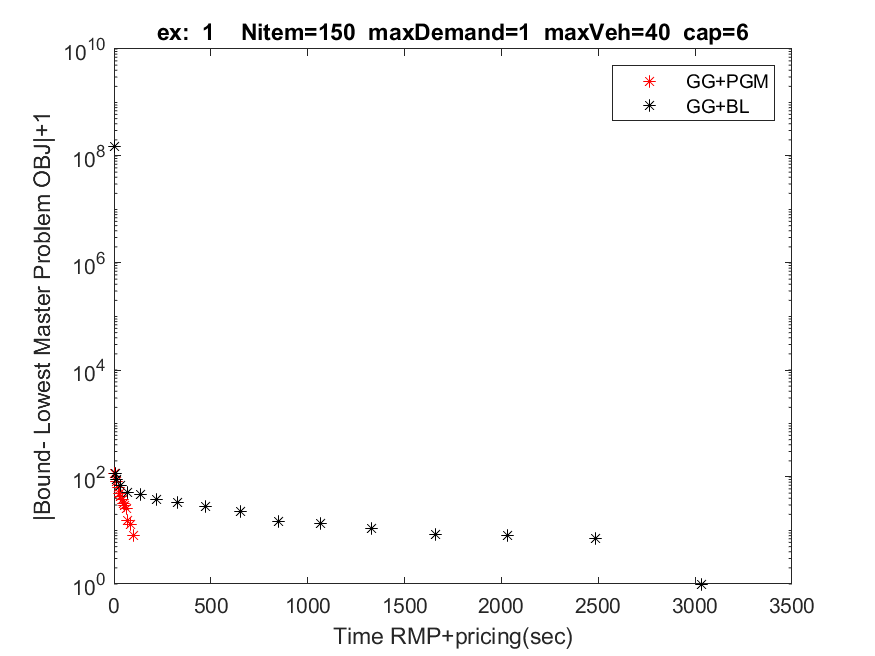}
	\includegraphics[width=0.49\linewidth]{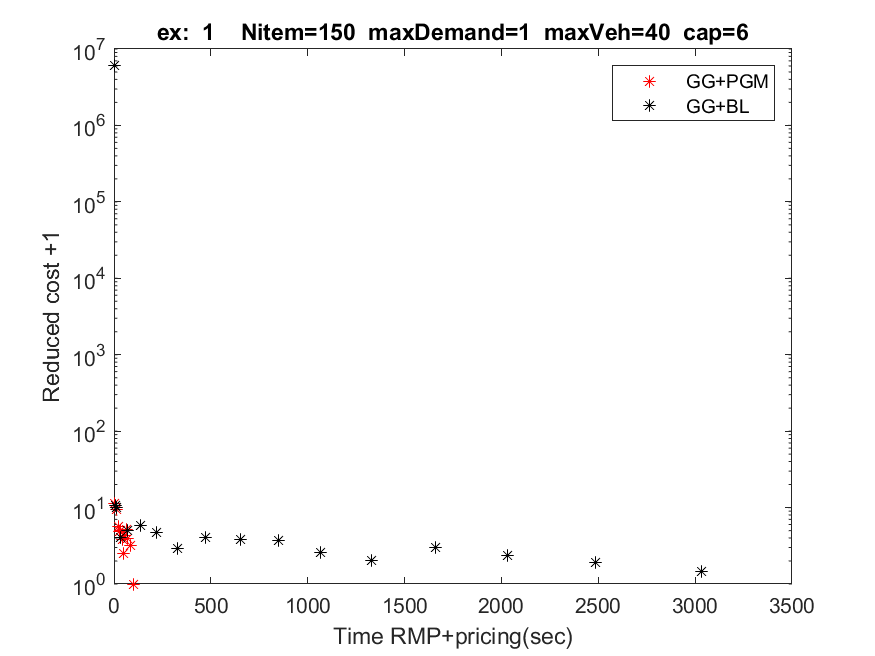}\\
		\includegraphics[width=0.49\linewidth]{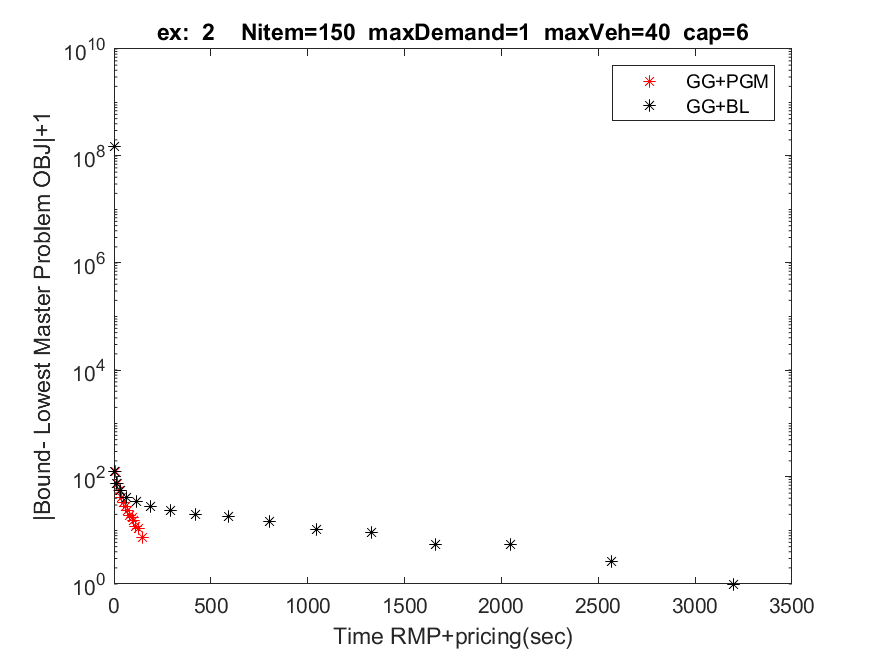}
	\includegraphics[width=0.49\linewidth]{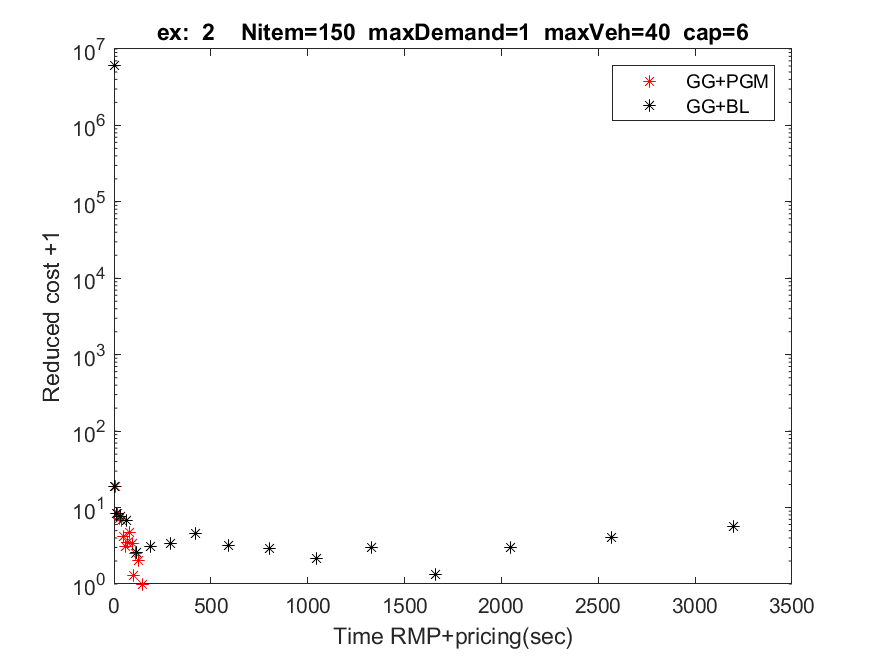}\\
	\includegraphics[width=0.49\linewidth]{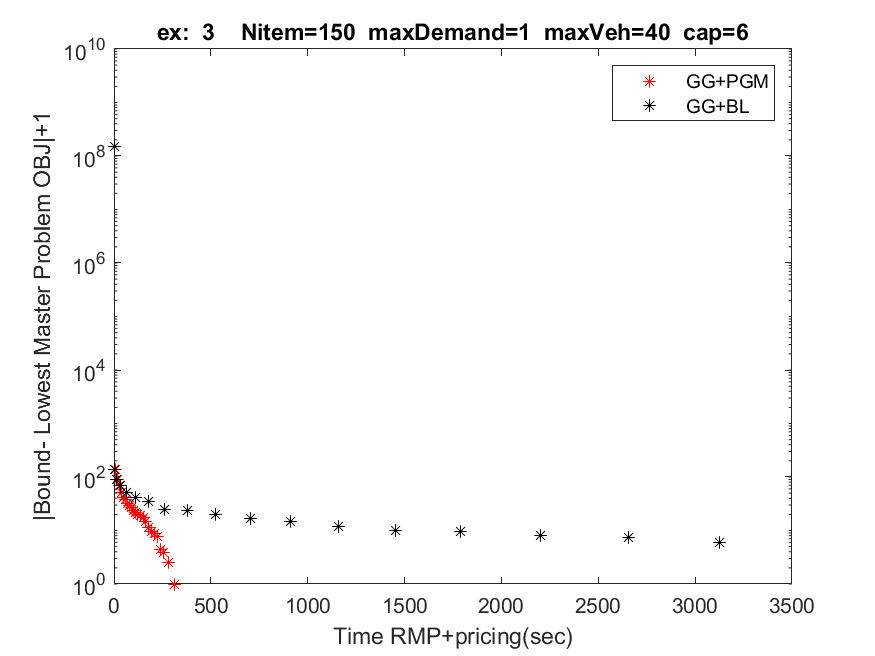}
	\includegraphics[width=0.49\linewidth]{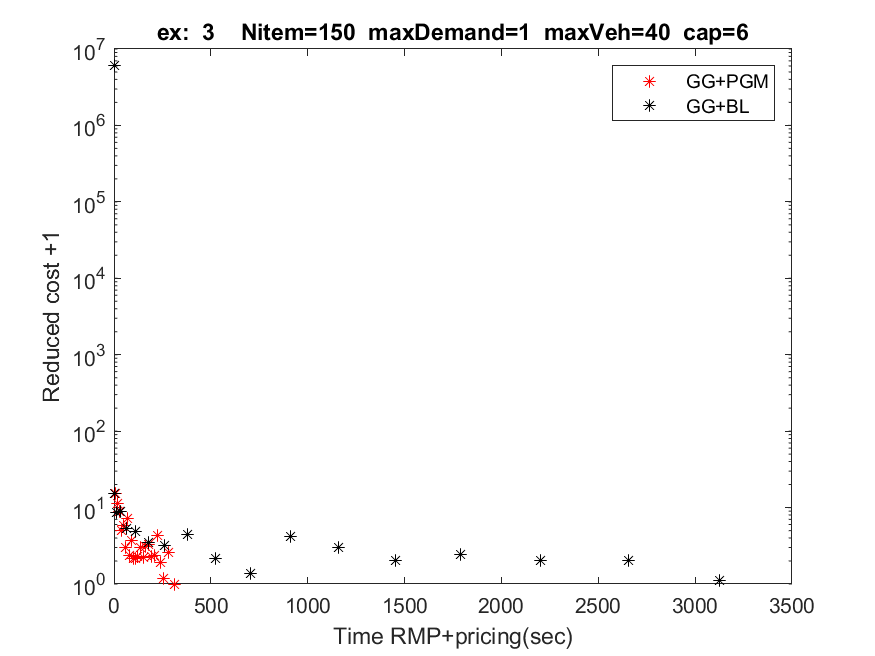}\\
	\caption{Results on individual problem instances as a function of time (sec) in semi-log scale. Individual dots show the a value (LP RMP or -reduced cost) on $y$ axis for the given algorithm at the iteration for that dot.  The left side provides the results RMP objective and the right side minus $1$ times the reduced cost. We add one to the $y$ values of all terms, which lets us use the semi-log scale. 
	}
	\label{fig:iter1}
\end{figure}

\begin{figure}[!hbtp]
	\includegraphics[width=0.49\linewidth]{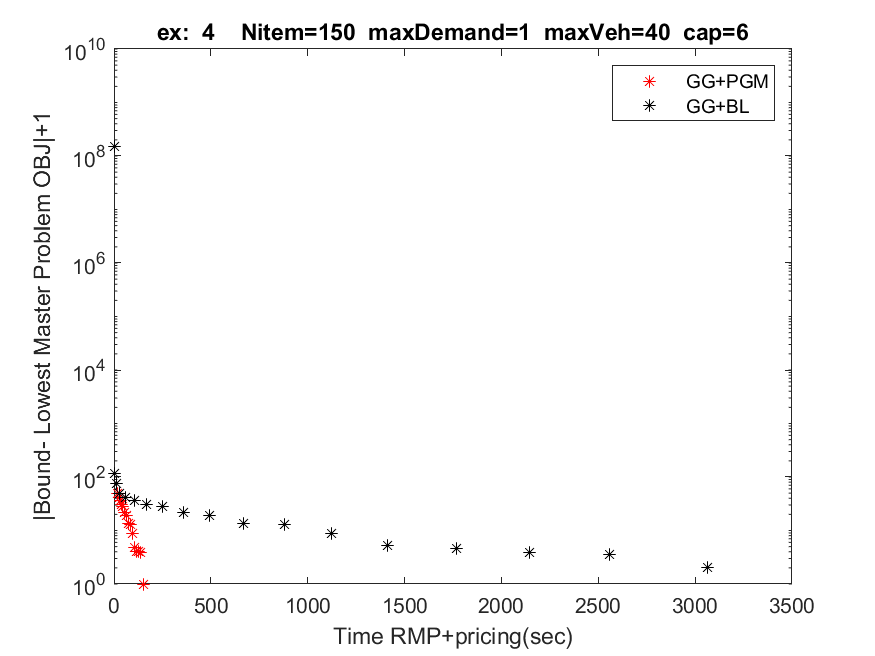}
	\includegraphics[width=0.49\linewidth]{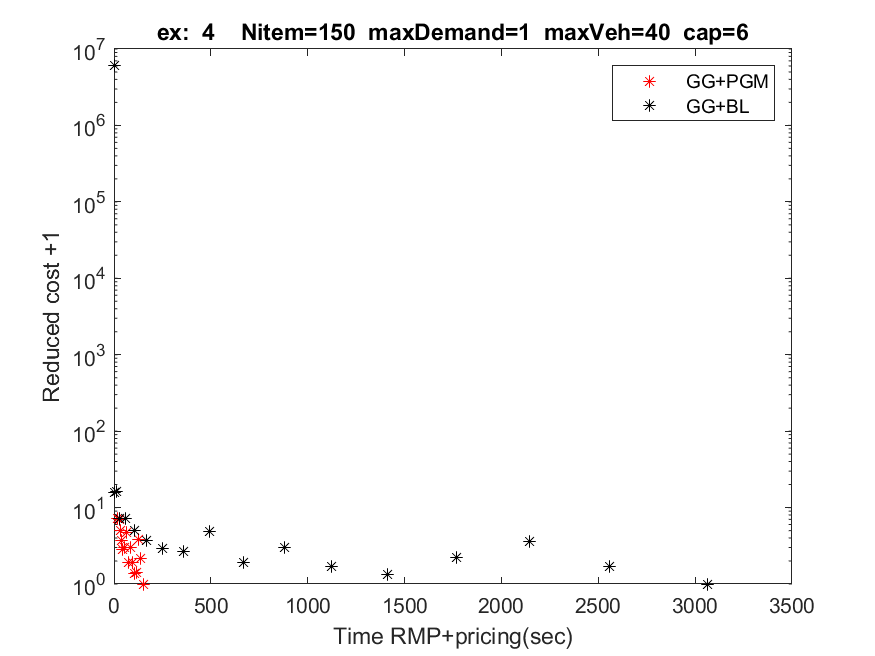}\\
		\includegraphics[width=0.49\linewidth]{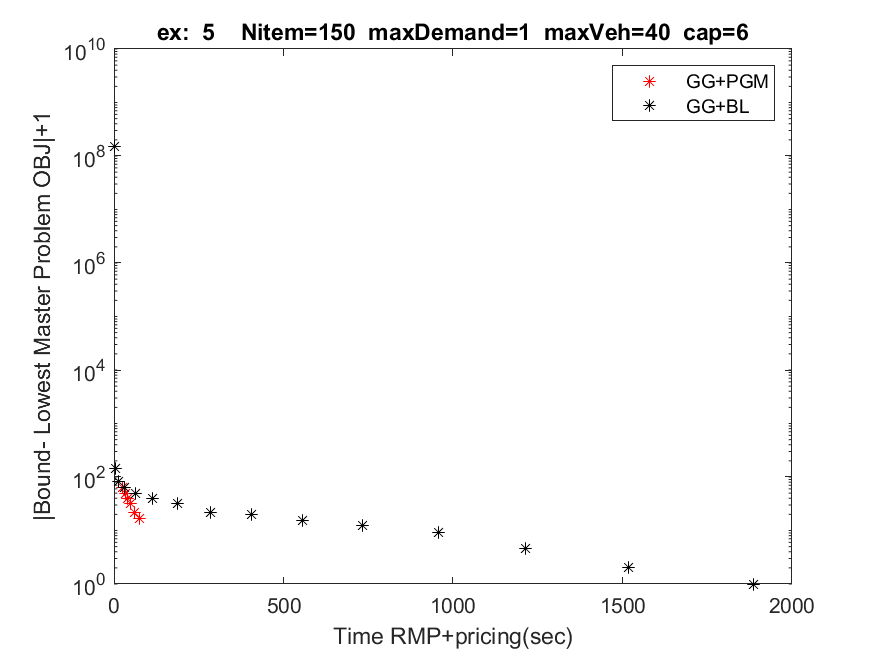}
	\includegraphics[width=0.49\linewidth]{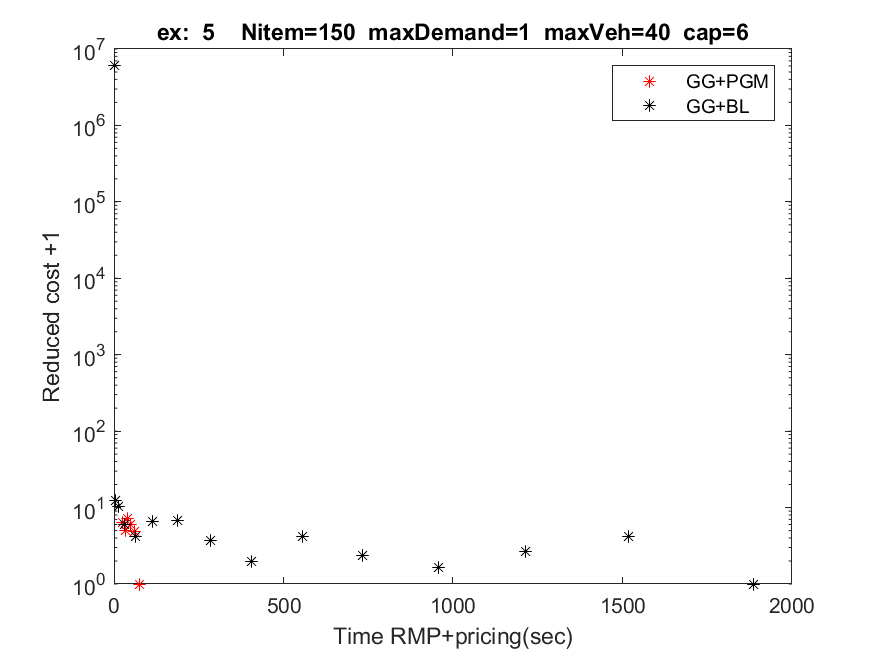}\\
	\includegraphics[width=0.49\linewidth]{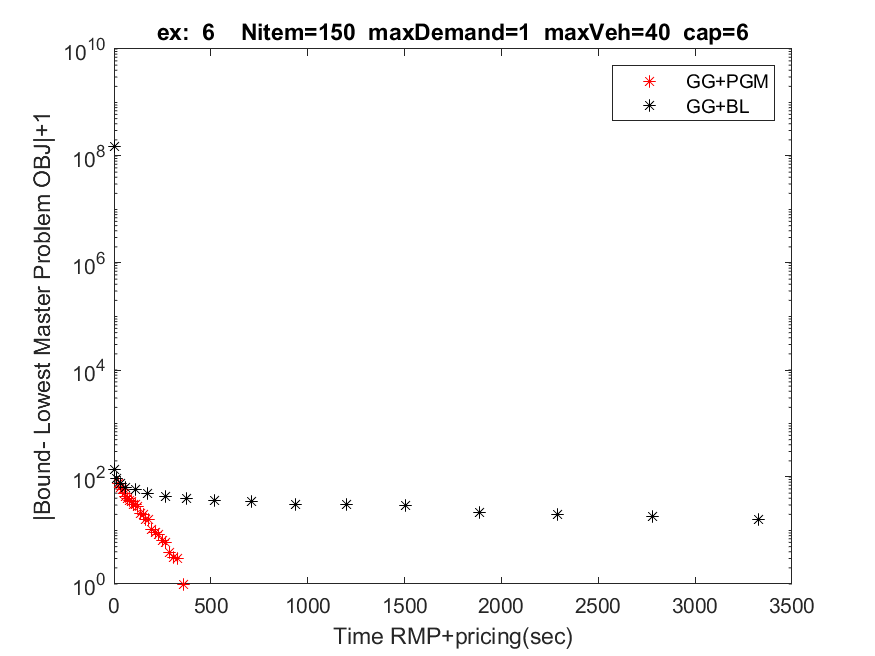}
	\includegraphics[width=0.49\linewidth]{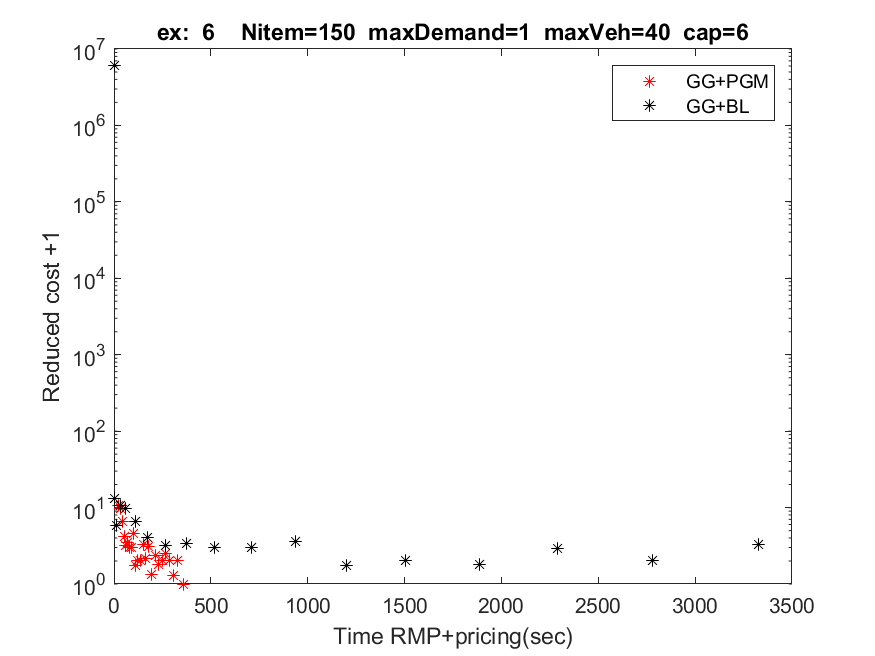}\\
	\caption{Results on individual problem instances as a function of time (sec) in semi-log scale. Individual dots show the a value (LP RMP or -reduced cost) on $y$ axis for the given algorithm at the iteration for that dot.  The left side provides the results RMP objective and the right side  minus $1$ times the reduced cost. We add one to the $y$ values of all terms, which lets us use the semi-log scale. 
	}
	\label{fig:iter2}
\end{figure}

\begin{figure}[!hbtp]
	\includegraphics[width=0.49\linewidth]{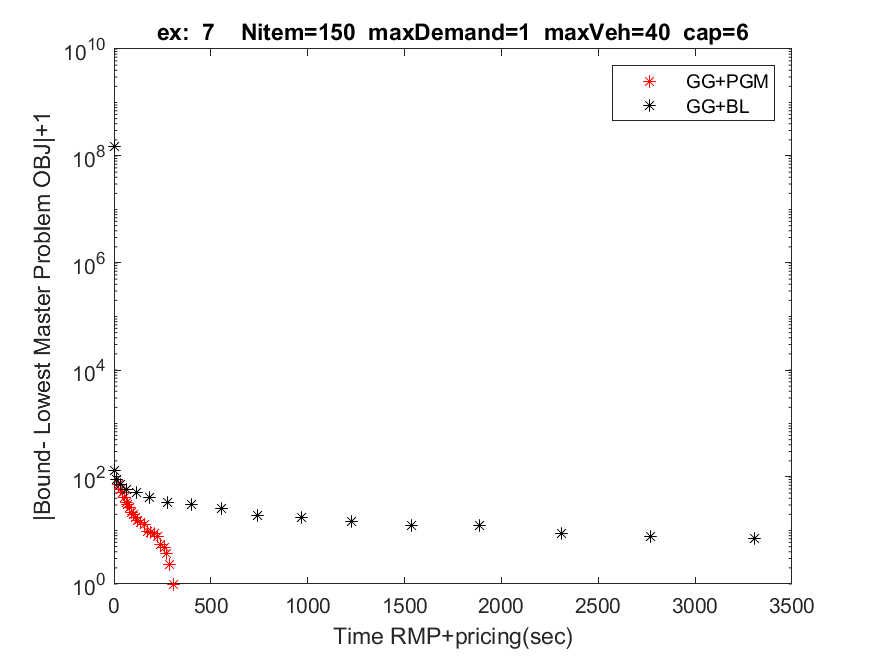}
	\includegraphics[width=0.49\linewidth]{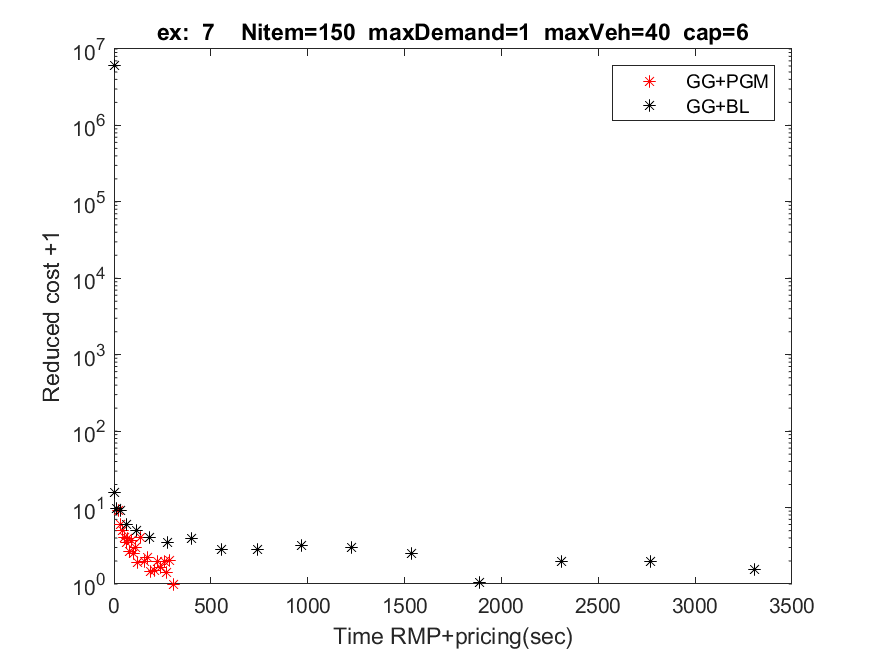}\\
		\includegraphics[width=0.49\linewidth]{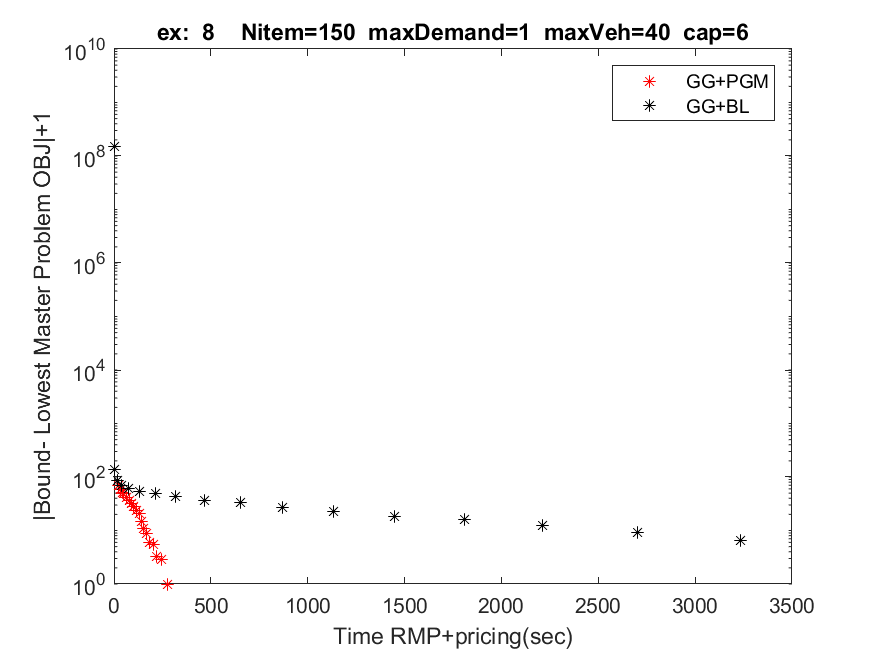}
	\includegraphics[width=0.49\linewidth]{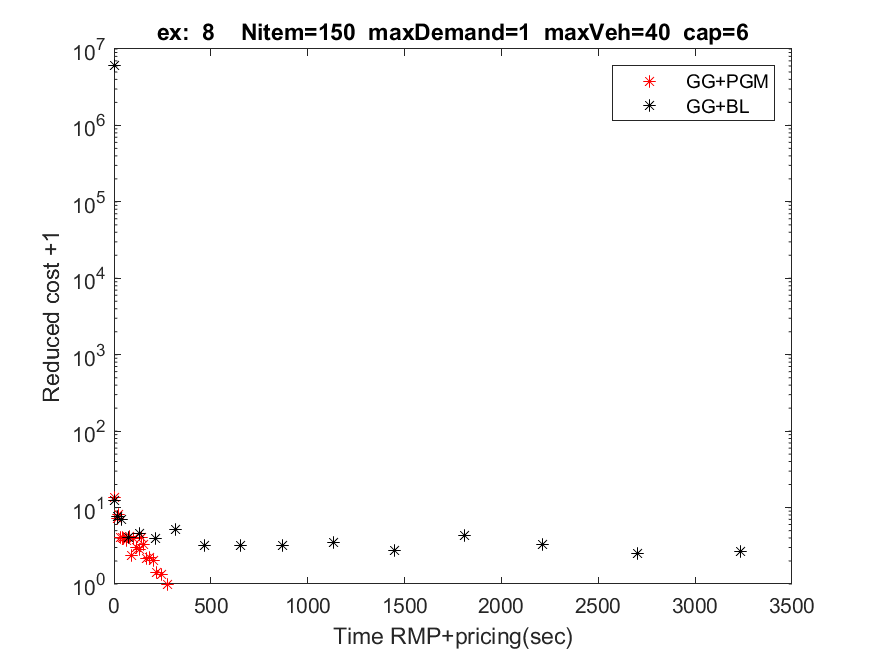}\\
	\includegraphics[width=0.49\linewidth]{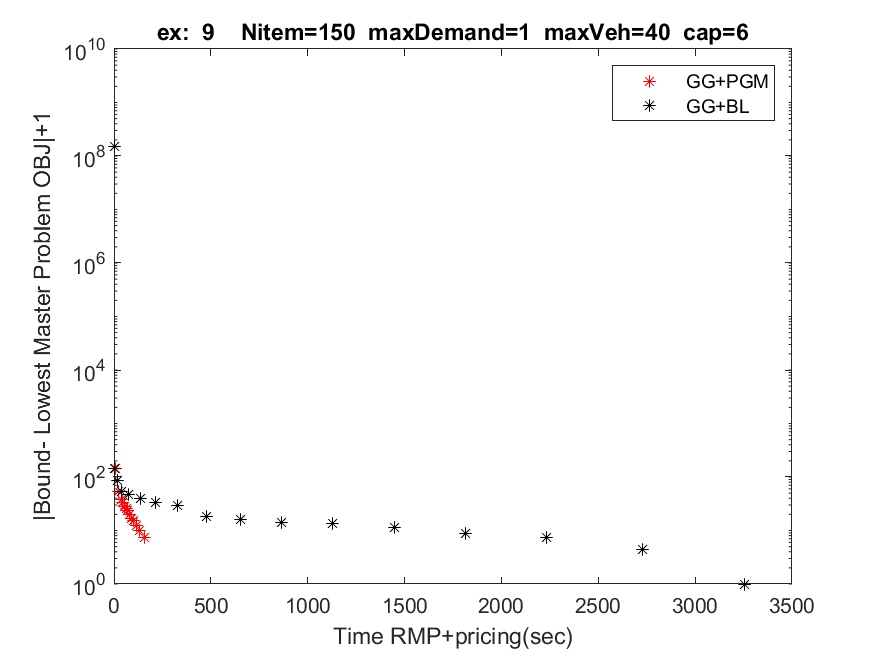}
	\includegraphics[width=0.49\linewidth]{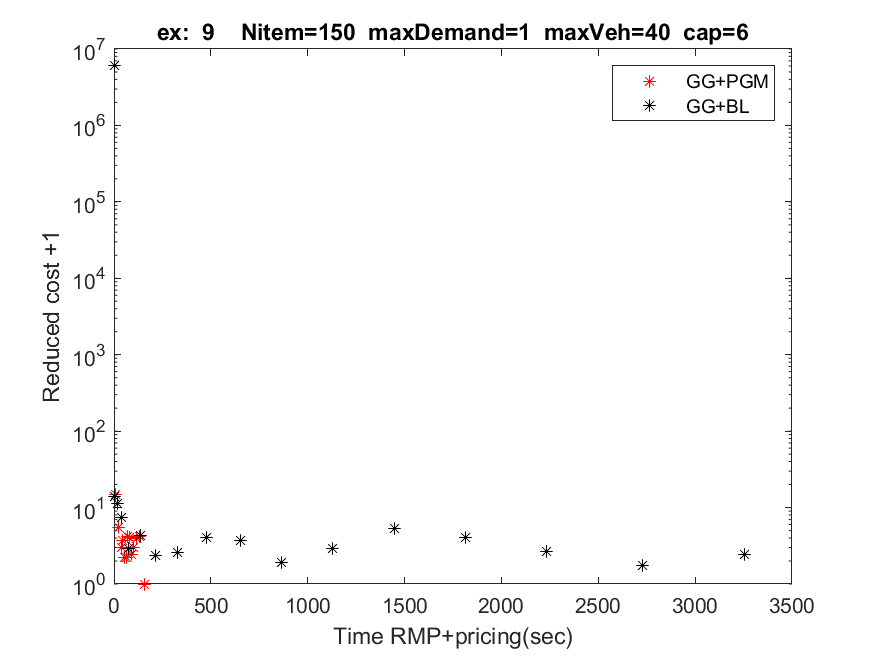}\\
	\caption{Results on individual problem instances as a function of time (sec) in semi-log scale. Individual dots show the a value (LP RMP or -reduced cost) on $y$ axis for the given algorithm at the iteration for that dot.  The left side provides the results RMP objective and the right side  minus $1$ times the reduced cost. We add one to the $y$ values of all terms, which lets us use the semi-log scale. 
	}
	\label{fig:iter3}
\end{figure}

We provide results on all our problem instances in Figs \ref{fig:iter1},\ref{fig:iter2},\ref{fig:iter3}.  On the left we provide the LP objective as a function of time and on the right the reduced cost  of the column generated in pricing as a function of time. We should note that the use of heuristic pricing not exact pricing can yield convergence to different LP values each time the problem instance is solved.  

\end{document}